\documentclass[regno,a4paper,11 pt]{amsart}
\usepackage[numbers]{natbib}
\usepackage[a4paper,top=4cm,bottom=2cm,left=2cm,right=2cm,bindingoffset=5mm]{geometry}

\usepackage[latin1]{inputenc}

\newtheorem{theorem}{Theorem}[section]

\newtheorem{proposition} [theorem]{Proposition}
\newtheorem{lemma}[theorem]{Lemma}
\newtheorem{corollary}[theorem]{Corollary}
\newtheorem{remark}[theorem]{Remark}

\newtheorem{Hypothesis}[theorem]{Hypotheses}
\makeatletter
\@addtoreset{equation}{section}

\usepackage[usenames]{color}
\definecolor{red}{rgb}{1.0,0.0,0.0}
\def\red#1{{\textcolor{red}{#1}}}

\newenvironment{system}
{\left\lbrace\begin{array}{@{}l@{}}}
{\end{array}\right.}

\newcommand{\tr}{\rm{tr}}

\def\ds{\begin{displaystyle}}
\def\eds{\end{displaystyle}}
\def\dis{\displaystyle }
\def\<{\langle }
\def\>{\rangle }

\newcommand{\R}{\mathbb{R}}
\newcommand{\E}{\mathbb{E}}

\newcommand{\e}{\varepsilon}

\newcommand{\mP}{\mathbb{P}}
\newcommand{\Fcal}{\mathcal{F}}     
\newcommand{\Lcal}{\mathcal{L}}
\newcommand{\calm}{\mathcal{M}}

\newcommand{\calu}{\mathcal{U}}
\newcommand{\calh}{\mathcal{H}}

\def\R{\mathbb R}

\def\E{\mathbb E}
\def\P{\mathbb P}

\def\F{\mathcal F}

\def\dim{\noindent \hbox{{\bf Proof.} }}

\title[SMP for SDEs with delay]{Stochastic maximum principle for equations with delay: going to infinite dimensions to solve the non-convex case}


\thanks{
The authors have been
supported by the Gruppo Nazionale per l'Analisi Matematica,
la Probabilit\`a e le loro Applicazioni (GNAMPA)
of the Istituto Nazionale di Alta Matematica (INdAM).}

\begin{document}

\author{Giuseppina Guatteri}
\address[G. Guatteri]{Dipartimento di Matematica, Politecnico di Milano. via Bonardi 9, 20133 Milano, Italia}
\email{giuseppina.guatteri@polimi.it}

\author{Federica Masiero}
\address[F. Masiero]{Dipartimento di Matematica e Applicazioni, Universit\`a di Milano-Bicocca. via Cozzi 55, 20125 Milano, Italia}
\email{federica.masiero@unimib.it}

\begin{abstract}
In this paper we develop necessary conditions for optimality, in the form of the stochastic Pontryagin maximum principle, for controlled equation with delay in the state
 and with control dependent noise, in the general case of controls $u \in U$ with $U$ not necessarily convex. The maximum principle is formulated by means of a first and second order adjoint BSDEs.
\end{abstract}

\maketitle



\smallskip \noindent


\section{Introduction}
In this paper we consider a controlled stochastic differential equation in $\R^d$ with delay in the state of the following form 
\begin{equation}\label{equazione_stato:ritardo-intro} 
\begin{system}
dx(t)=b(t, x(t),x_t, u(t))dt + \sigma(t,x(t),x_t,u(t))dW(t), \\
x(\theta)= x_0(\theta), \quad u(\theta)=\eta(\theta),\qquad \theta \in [-d, 0].
\end{system}
\end{equation}
Here and throughout the paper we use the notation  $x_t(\theta) = x(t+\theta),\,  \theta \in [-d,0]$ to denote the past trajectory of $x$ from $t-d$ up to time $t$.  
We consider admissible controls $u$, that are progressively measurable and square integrable processes taking values in a non convex set $ U\subset \R^k$: in this case the stochastic maximum principle can be formulated in terms of first and second order adjoint equations.
\newline We are able to allow a quite general dependence on the past trajectories $x_t$  of the state, namely the drift  and diffusion can be written as 
\begin{align}\label{barf-barg-intro}
 b(t,x,y,u)=b (t,x, \int_{-\mathtt{d}}^0 y(\theta )\mu_b(d\theta), u(t )), 
 \quad \sigma(t,x,y,u)=\sigma (t, x(t),\int_{-\mathtt{d}}^0 y(\theta )\mu_\sigma(d\theta), u(t) ),
 \end{align}
 where $\mu_b,\mu_\sigma$ are finite regular measures on $[-d,0]$, and $b(t,\cdot,\cdot,\cdot)$ and  $\sigma(t,\cdot,\cdot,\cdot)$ are Lipschitz continuous and twice differentiable.
Associated to equation (\ref{equazione_stato:ritardo-intro}) we consider the cost functional  
\begin{equation}\label{costo_intro}
J(u(\cdot)) = \E \int_0^T \ell(t,x(t),x_t, u(t)) dt + \E\, h(x(T),x_T).
\end{equation}
that we have to minimize over all admissible controls. Also in the running cost $\ell$ and in the final cost $h$ we can allow a general dependence on the past trajectories $x_t$ and $x_T$: there exist finite regular measures $\mu_\ell,\,\mu_h$ such that the current cost and the final cost can be written as follows \vspace{-0.2truecm}
\begin{align}\label{barl-barh-intro}
\ell(t,x,y, u)= \ell (t,x,\int_{-\mathtt{d}}^0 y(\theta )\mu_\ell(d\theta),u),\quad 
 h(x,y)=h (x,\int_{-\mathtt{d}}^0 y(\theta )\mu_h(d\theta)), \; \forall \,y \in \, C([-d,0], \R^d),
\end{align}
We first consider the case when the measures $\mu_b,\mu_\sigma, \mu_\ell, \mu _h$ are absolutely continuous with respect to the Lebesgue measure, and then, by an approximation procedure we are able to consider general regular measures, among others the Dirac measures, so that in equation \eqref{equazione_stato:ritardo-intro} we are able to consider also the relevant case of pointwise delay.
 Notice that in the state equation \eqref{equazione_stato:ritardo-intro} and in the cost functional \eqref{costo_intro} we do not consider delay in the control (see respectively \eqref{barf-barg-intro} and \eqref{barl-barh-intro}): this choice has been done for the sake of simplicity, and we stress the fact that with our techniques we can cover also the case of delay in the control, as done in \cite{GuaMas}.

In the case of control processes taking values in a convex set $U$, the problem has already been treated: for the case of pointwise delay see e.g. \cite{ChenWu}, see also \cite{ChenWuYu} where the linear quadratic case is considered and 
\cite{ChenWu2011} where the authors prove also existence of optimal controls; in the case of more general dependence on the past the problem has been treated in \cite{GuaMas}, see also \cite{GMO} where SPDEs with delay are studied, and in \cite{HuPe}, where the authors
formulate the stochastic maximum principle by functional analysis methods.
\newline In all the aforementioned papers, except the last one,  the stochastic maximum principle for problems with delay is formulated by means of anticipated BSDEs, that have been introduced in \cite{PengYang} and generalized in \cite{YangEll}. 

In the present paper, we consider the general case of control dependent noise when the control processes take values in a set $U $ not necessarily convex and we apply the so called spike variation method, consisting in perturbing the optimal control only in a measurable set $E_\e\subset [0,T]$ and considering the variation of the state with respect to this perturbation
of the control. 

\noindent The spike variation method has been introduced in \cite{Peng} for a state equation without delay, see also \cite{YongZhou}.
\noindent In \cite{Peng} the maximum principle has been formulated 
by means of two backward stochastic differential equations (BSDEs in the following): towards the case of $U$ convex, besides the (first order) adjoint BSDE, which is the adjoint equation of the variation of the state, the authors introduce a second order adjoint BSDE,  which is the adjoint equation of the equation for the square of variation of the state.
\newline In the recent papers \cite{MeShi} and \cite{MeShiWa} the authors consider the case with delay state equation, with control dependent noise and with controls taking values in a not necessarily convex set: in \cite{MeShi} the authors solve the problem by imposing very strong conditions, namely they introduce an additional BSDE requiring its solution to be zero; in \cite{MeShiWa} the authors are able to treat a much more general case considering backward stochastic Volterra equations with two kind of dependence on the past trajectory: pointwise delay and integral dependence on the past with exponential density. 

Here we propose a direct method based on the infinite dimensional reformulation of the problem, in the Hilbert space $H=\mathbb{R}^d\oplus L^{2}\left([-\mathtt{d},0] ,\mathbb{R}^d\right)$, where the first component is the present state and the second component is its past trajectory, see Section \ref{sec-reformul} for more details. In such a way the problem is reformulated as an infinite dimensional control problem without delay, we can so proceed similarly to \cite{FuHuTe-CR} and \cite{FuHuTe-AMO}, where the authors prove an infinite dimensional version of the stochastic maximum principle. In the case when the dependence on the past is given by measures with a density with respect to the Lebesgue measure this approach completely solves the problem. 
\newline In the case when the dependence on the past is given by general measures, in order to work with an infinite dimensional reformulation of the problem, we should work in the Banach space $E=\mathbb{R}^d\oplus C\left([-\mathtt{d},0] ,\mathbb{R}^d\right)$, but this would bring together a first order adjoint equation for a pair of processes taking values in $E^*=\mathbb{R}^d\oplus \calm \left([-\mathtt{d},0]\right)$, where $\calm\left([-\mathtt{d},0]\right)$ is the space of regular measures on $[-d,0]$, and a second order adjoint equation with values in a subspace of $L(E^*, E^*)$. Since working in this spaces would be very complicated, we proceed as follows: at first we do not perform any infinite dimensional reformulation, we consider a first order adjoint equation  which turns out to be an  anticipated BSDE; in order to treat the terms containing the square of the present state and of the past trajectory we approximate the first and  second order variation of the state by approximating the general regular measure by means of a measure with density, see \cite{GMO} where a similar approximation of measures is performed. At this point we can write an approximating second order adjoint equation for the pair of processes $(P^n,Q^n)$ with values in $S^2(H)$, the Hilbert space of Hilbert-Schmidt symmetric operators from $H$ to $H$. Notice that we are not able to prove the convergence of $(P^n,Q^n)$ to the solution of a BSDE with values in $S^2(H)$ nor we are able to prove the convergence of the pair $(P^n,Q^n)$. It is crucial that, in order to formulate the stochastic maximum principle,  what matters in only the finite dimensional component $P^n_{0,0}$: we are able to prove that $P^n_{0,0}$ converges, and we are able to characterize its limit $P_{0,0}$, see Proposition \ref{Conv P00} and Theorem \ref{teoP00}.

The paper is organized as follows. In section 2 we set notation and assumptions, section 3 is devoted to the reformulation of the problem in a suitable infinite dimensional space, to the statements on the first and second adjoint equations and to the formulation of the maximum principle when the measures are absolutely continuous with respect to the Lebesgue measure. Section 4 is devoted to the case of general measures. In the Appendix we discuss an application and we give some technical proofs.

\section{Assumptions and preliminaries}\label{sec:AssPrel}

Let $(\Omega, \Fcal, \P)$ be a complete probability space
and let $W(t)$ be a  $m$ dimensional brownian motion. We endow $(\Omega, \Fcal, \P)$ with the natural filtration $(\mathcal{F}_t)_{t\geq 0}$ generated by $W$ and
augmented in the usual way with the family of $\mP$-null sets of $\mathcal{F}$. A finite time horizon $T>0$ is fixed, we denote by $H$ a real and separable Hilbert space endowed with a norm denoted by $|\cdot|_H$ ($|\cdot|$ if no confusion is possible) induced by the inner product  $\langle \cdot, \cdot \rangle_H$ ( $\langle \cdot, \cdot \rangle$ again if no confusion is possible.

Let now $E, \vert \cdot \vert $ be a Banach space; we introduce the following spaces for any $p \geq 1$ and for every $ 0\leq r\leq s\leq T$:
\begin{itemize}
\item $L^p(\Omega, \F_s;E)$ the set of $\F_s$ measurable random vectors in $H$;
\item $L^p_{\mathcal{F}}(\Omega\times [r,s];E)$, the set of all $(\mathcal{F}_t)$-progressively measurable processes with values in $E$ such that\vspace{-0.2truecm}
\[ \|X\|_{L_{\mathcal{F}}^p(\Omega\times[r,s];E)} = \left( \E \int_r^s |X(t)|^p dt\right)^{1/p} < \infty; \]
\item $L^p_{\mathcal{F}}(\Omega; C([r,s];E))$, the set of all $(\mathcal{F}_t)$-progressively measurable processes with values in $E$ such that\vspace{-0.2truecm}
\[ \|X\|_{L^p_{\mathcal{F}}(\Omega; C([r,s];E
))} = \Big( \E \sup_{t \in [r,s]} |X(t)|^p \Big)^{1/p} < \infty. \]
\end{itemize}

\subsection{Formulation of the problem}\label{sec-density}
Let $U \subset \R^k$ be a non empty set.
By admissible control we mean any  process  $u: \Omega\times [0,T] \to U $ that is $(\mathcal{F}_t)$-progressive measurable
and we denote by $\calu_{ad} $  the set of admissible controls.
\newline As already noticed, we do not make any convexity assumptions on the space $U$: this allows e.g. to consider the case of $U= \{0,1\}$, or more generally
of $U$ any discrete space.

We study the following stochastic state equation with delay\vspace{-0.2truecm}
\begin{equation}\label{equazione_stato} 
\begin{system}
dx(t)=  b\big(t,x(t), \dis\int_{-\mathtt{d}}^0 x(t+\theta )\mu_b(d\theta), u(t )\big)+\sigma \big(t, x(t),\dis\int_{-\mathtt{d}}^0 x(t+\theta )\mu_\sigma(d\theta), u (t)\big)dW(t), \\
x(0)= x_0; \,\qquad x(\theta)= x_1(\theta)\qquad \theta \in [-d, 0);
\end{system}
\end{equation}
where $x(t),\, x_0 \in \R^d$, while $x_1 \in C([- \mathtt{d},0]; \R^d)$. 
 
We make the following assumptions on the coefficients:
\begin {Hypothesis}\label{ipotesi} 
 \begin{itemize}
\item[]
\medskip
\item [(i)] The functions $b$ and $\sigma$:\vspace{-0.2truecm}
\[  b: [0,T]  \times \R^d \times \R^d \times \R^k \to \R^d, \qquad  \sigma: [0,T]  \times \R^d \times \R^d \times \R^k \to \R^{d \times m}\]
are continuously differentiable with  bounded  derivatives with respect to $(x,y)$  up to the second order 
 and there exist positive constants $M$ and $L$, such that for $\psi = b,\sigma$:
\begin{equation}\label{Lip}
|\psi(t,x,y,u) - \psi(t,x,y,u')| \leq L|u-u'|, \\ \forall t \in [0,T], u,u' \in U.
\end{equation}
\begin{equation}\label{Bound}
|\psi(t,0,0,u)| \leq  M \qquad \forall t \in [0,T], u \in U.
\end{equation}
\item[(ii)] The measures $\mu_b,\, \mu_\sigma$ are finite regular measures on $[-\mathtt{d},0]$. 
\end{itemize}
\end{Hypothesis}
 \medskip
 
To the controlled state equation with delay \eqref{equazione_stato}, we associate the following cost functional 
\begin{equation}\label{costo}
J(u) :=\E\int_0^T \ell(t,x(t), \dis\int_{-\mathtt{d}}^0 x(t+\theta )\mu_\ell(d\theta)),u(t)) \, dt + \E(h(x(T)), \dis\int_{-\mathtt{d}}^0 x(T+\theta )\mu_h(d\theta))
\end{equation}
where we assume:
\begin{Hypothesis}\label{ipotesi-costo}
\begin{itemize}\item[]\medskip
\item [(i)]The functions $\ell(t,x,y,u)$ and $h(x,y)$:
\[  \ell: [0,T]\times \R^d \times \R^d \times \R^k \to \R, \qquad h: \R^d \times \R^d  \to \R\]
are continuously differentiable with bounded derivatives with respect to $(x,y)$
 and there exist positive constants $M$ and $L$, such that:
\begin{equation}\label{Lipcosto}
|\ell(t,x,y,u) - \ell(t,x,y,u')| \leq L|u-u'|, \\ \forall t  \ \in [0,T], \ x,y \in \R^d \  u,u' \in U.
\end{equation}
\begin{equation}\label{Boundcosto}
|\ell(t,0,0,u)| \leq  M \qquad \forall t \in [0,T], u \in U.
\end{equation}
\item [(ii)]The measures $\mu_\ell,\, \mu_h$ are finite regular measures on $[-\mathtt{d},0]$. 
\end{itemize}
\end{Hypothesis}
\medskip
\noindent The control problem consists in minimizing over all admissible controls the cost functional given in (\ref{costo}).

Next we recall  existence and uniqueness results for equation \eqref{equazione_stato}, see e.g. \cite{Moha}
\begin{theorem}\label{teo-mohammed}
UnderHypothesis \ref{ipotesi}, $\forall x_1 \in  C([-\mathtt{d},0]; \R^d)$ and $ \forall \,u \in \mathcal{U}$ one has that there exists a unique solution $x$
of \eqref{equazione_stato} with $x_0= x_1(0)$ and there exists a constant $c>0$ such that:\vspace{-0.2truecm}
\begin{equation} \label{stima_dati_stato}
\E \sup_{t \in [-\mathtt{d},T]} |x(t)|^2 \leq c (1+ |x_1|^2_{C([-\mathtt{d},0]; \R^d)})
\end{equation}
\end{theorem}

\section{The case of past dependence through measures with a density: infinite dimensional reformulation}
\label{sec-reformul}

We first consider the case of finite regular measures $\mu_b,\,\mu_\sigma,$ in the state equation, absolutely continuous with respect to the Lebesgue measure, namely we assume that:
\begin{Hypothesis}\label{ip-densita-eqstato}
The finite regular measures $\mu_b$ and $\mu_\sigma,$ are absolutely continuous with respect to the Lebesgue measure with densities square integrable in $[-\mathtt{d},0]$ which are denoted respectively by $f_b$ and $f_\sigma$.
\end{Hypothesis}
\noindent Notice that under Hypothesis \ref{ip-densita-eqstato}, as initial past trajectory in equation \eqref{equazione_stato} we could handle also $x_1\in L^2([-\mathtt{d},0], \R^d)$. 

We are going to reformulate equation \eqref{equazione_stato} as an infinite dimensional stochastic evolution equation without delay, indeed we notice that due to the dependence at time $t$ on the past trajectory $(x (t+\theta))_{\theta\in[-\mathtt{d},0]}$, the problem is intrinsic infinite dimensional: it is reformulated in an infinite dimensional space where both the present and the past trajectory evolve. We consider the classical reformulation in the Hilbert space of square integrable trajectories, see e.g. \cite{CM}, \cite{DM} and \cite{W}, and for the stochastic case with multiplicative noise\cite{BGPZ} and\cite{BGP}.
\newline Coming into the details of the reformulation, we define the Hilbert space $H=\mathbb{R}^d\oplus L^{2}\left([-\mathtt{d},0] ,\mathbb{R}^d\right)$. 
It turns out that for $t\geq 0 $\vspace{-0.2truecm}
\begin{equation}\label{semigr-delay}
e^{t A}:H\longrightarrow H,\quad e^{t A}\left(
\begin{array}[c]{@{}c@{}}%
x_{0}\\
x_{1}%
\end{array}
\right) =\left(
\begin{array}[c]{@{}c@{}}%
x_0 \\
x_01_{[-t,0]}(\cdot)+ x_1(\cdot+t)1\cdot_{[-d,-t ] }%
\end{array}
\right),
\end{equation}
defines a $C_{0}$-semigroup in $H$; see, e.g.,~\cite{DM} and \cite{W}.
The infinitesimal generator $A$ of $e^{t A},t\geq 0$, is
given by\vspace{-0.1truecm}
\begin{equation}\label{A}
\begin{array}[c]{l}%
\mathcal{D}\left( A\right) =\left\{ \left(
\begin{array}[c]{@{}c@{}}%
x_{0}\\

x_{1}%
\end{array}
\right) \in H,x_{1}\in H^{1}\left( \left[ -d,0\right] ,\mathbb{R}%
^{d}\right) ,x_{1}\left( 0\right) =x_{0}\right\},\quad
Ah=A\left(
\begin{array}[c]{@{}c@{}}%
x_{0}\\

x_{1}%
\end{array}
\right) =\left(
\begin{array}[c]{@{}c@{}}%
0 \\

dx_{1}/d\theta
\end{array}
\right) .
\end{array}
\end{equation}
By setting\vspace{-0.2truecm}
\begin{equation}\label{G}
G:\mathbb{R}^{n}\longrightarrow H,\text{ }G=\left(
\begin{array}[c]{@{}c@{}}%
I\\
0
\end{array}
\right) ,\text{ }X(t)=\left(
\begin{array}[c]{@{}c@{}}%
x\left( t\right) \\
(x(t+ \theta))_{ \theta \in [-\mathtt{d},0]})%
\end{array}
\right),\quad \text{and}\quad X_{0}=\left(
\begin{array}[c]{@{}c@{}}%
x_{0}\\
x_{1}%
\end{array}
\right),
\end{equation}
problem (\ref{equazione_stato}) can be rewritten as an abstract evolution equation in $H$ 
\begin{equation}
\left\lbrace \begin{array}
[c]{l}
dX(t)  =AX(t) dt+GB(t,X(t),u(t))\, dt +G\Sigma(t,X(t),u(t))\,dW_t
,\text{ \ \ \ }t\in[ 0,T] \\
X(0)  =x=\left(\begin{array}{l}x_0\\x_1\end{array}\right).
\end{array}
\right.   \label{eq-astr}%
\end{equation}
where $B:[0,T]\times H\times \R^k\rightarrow\R^d,\;\Sigma:[0,T]\times H\times \R^k\rightarrow\R^{d\times m} $, and for $t\in[0,T],\,Y=\left(\begin{array}{l}
Y_0\\
Y_1
\end{array}\right)\in H,\,u\in R^k$, $B$ and $\Sigma$ are given by\vspace{-0.2truecm}
\begin{equation}\label{B-Sigma}
B(t,Y,u):=b(t,Y_0, \dis\int_{-\mathtt{d}}^0Y_1(\theta)f_b(\theta)\,d\theta),u),\;
\Sigma(t,Y,u):=\sigma(t,Y_0, \dis\int_{-\mathtt{d}}^0Y_1(\theta)f_\sigma(\theta)\,d\theta),u)
\end{equation}
Taking the integral mild form of (\ref{eq-astr}) we have
\begin{equation}
X(t)  =e^{tA}x+\int_0^te^{(t-s)A}GB(s,X(s),u(s))\,ds+\int_0^te^{(t-s)A}G\Sigma(s,X(s),u(s))\,dW_s
,\text{ \ \ }t\in[ 0,T]. \\
  \label{eq-astr-mild}%
\end{equation}
In analogy toHypothesis \ref{ip-densita-eqstato}, we assume that the measures $\mu_\ell,\,\mu_h$ are absolutely continuous with respect to the Lebesgue measure, namely we assume that
\begin{Hypothesis}\label{ip-densita-costo}
The finite regular measures $\mu_\ell$ and $\mu_h$ are absolutely continuous with respect to the Lebesgue measure with densities square integrable in $[-\mathtt{d},0]$ which are denoted respectively by $f_\ell$ and $f_h$.
\end{Hypothesis}
\noindent Under this assumption, also the cost functional can be reformulated in an abstract way as\vspace{-0.2truecm}
\begin{equation}\label{costo-abstract}
J(u) :=\E\int_0^T L(t,X(t),u(t)) \, dt + \E(H(X(T)) )
\end{equation}
where for $Y=\left(\begin{array}{l}
Y_0\\
Y_1
\end{array}\right)\in H$, $L$ and $H$ are given by\vspace{-0.2truecm}
\begin{equation}\label{L-H}
L(t,Y,u):=\ell(t,Y_0, \dis\int_{-\mathtt{d}}^0 Y_1(\theta )f_\ell(\theta )d\theta),u),  \quad  H(Y)=h(Y_0, \dis\int_{-\mathtt{d}}^0 Y_1(\theta )f_h(\theta) d\theta)
\end{equation}
So the problem of finding necessary conditions for the existence of optimal controls when the space of controls is not convex and the state equation is a stochastic delay equations, has turned into a problem related to an infinite dimensional stochastic evolution equations. In \cite{FuHuTe-CR} and \cite{FuHuTe-AMO} the authors study an analogous general stochastic maximum principle for stochastic partial differential equations, where $A$ is an analytic operator: we follow the line of those papers, even if in the present paper $A$ turns out to be the generator of a strongly continuous semigroup of linear operators and we do not assume any regularizing properties on the semigroup $(e^{tA})_{t>0}$.
\subsection{First and second order variation of the state.}
Let us assume that ${u}$ is an optimal control and ${X}$ the corresponding optimal trajectory, we introduce the {\em spike} variation:
let $\e >0$ and $E_\e \subset [0,T]$, such that $\lambda(E_\e)= \e$, where $\lambda$ is the Lebesgue measure.
Let $v \in U$ 
we set\vspace{-0.2truecm}
\begin{equation}\label{spike_def}
u^\e(t) = \left\{ \begin{array}{ll}
{u} (t)  & t \in [0,T] \setminus E_\e \\
v & t \in E_\e
\end{array}\right.
\end{equation}
We will denote by $X^\e$ the solution of \eqref{eq-astr} corresponding to $u^\e$.
For semplicity we  use the following notation \vspace{-0.5cm}
\begin{align*}
&\qquad \qquad\phi(t):=\phi(t,X(t),u(t)), \qquad
\phi^\e(t):= \phi(t,X^\e(t),u^\e(t)), \\\vspace{1mm}
& \nabla_\rho\phi:=\phi_\rho(t)= \phi_\rho(t,X(t),u(t)),\qquad\delta \phi(t):= \phi(t,X(t),u^\e(t))- \phi(t,X(t),u(t)),
\end{align*}
where $\phi= B ,\,\Sigma ,\, L, H$ and $\rho= X= (X_0,X_1)$.
We introduce the $H$ valued equation for the first variation of the state 
\begin{equation}\label{equazione_varI} \left\{
\begin{array}{l}
dY^\e(t)  = AY^\e(t)\,dt  +GB_X(t) Y^\e(t)\,dt  + \sum_{i=1}^m G \Sigma_X^i(t) Y^\e (t) d\, W ^i(t)  + \sum_{i=1}^m G \delta \Sigma ^i(t)d\, W ^i(t),  \; t \in [0,T],\\  \\
Y(0) =0\in H. 
\end{array}
\right.
\end{equation}
We introduce also the $H$ valued equation for the second variation
\begin{equation} \label{equazione_varII}\left\{
\begin{array}{rl}
dZ^\e(t) &  = [AZ^\e(t) +G  B_X(t) Z^\e(t)  +    G\delta B (t)  + \frac{1}{2} G B_{XX}(t) Y^\e(t)^2]\, dt  \\  \\& + 
 \sum_{i=1}^m [ G\frac{1}{2} \Sigma_{XX} ^i(t)  Y^\e(t)^2 + G\Sigma^i_X(t) Z^\e (t)  + 
G\delta  \Sigma^i_{X}(t) Y^\e(t)]d\, W^i(t),  \quad t \in [0,T], \vspace{0.3cm} \\ 
 Z^\e_0  & =\,0 \in H.
\end{array} \right.
\end{equation}
 where  for every $h = (h_0, h_1(\cdot)), k =  (k_0, k_1(\cdot)) \in H$ we have that \vspace{-0.2truecm}
 \begin{align*}
 B_{X}(t) h  &=\left( b ^ j_x(t) h_0 +  \int_{-\mathtt{d}}^0 b ^j_y(t)f_b(\theta)  h_1(\theta)\, d \theta\right)_{1 \leq j \leq d} \\ 
   \Sigma^i_X(t) h & = \left(  \sigma^{j,i}_x(t) h_0 +  \int_{-\mathtt{d}}^0 \sigma^{j,i}_y(t) f_\sigma(\theta)  h_1(\theta)\, d \theta \right)_{1 \leq j \leq d} =  \Big( \Sigma_{X} ^{j,i}(t) h \Big)_{1 \leq j \leq d }, \qquad   i: 1, \dots  m,\\
\delta  \Sigma^i_X(t) h &=   (\sigma^{j,i}_x(t) - \sigma^{j,i, \e}_x(t))  h_0 + \int_{-\mathtt{d}}^0( \sigma^{j,i}_y(t)  -  \sigma^{j,i, \e}_y(t))f_\sigma(\theta)  h_1(\theta)\, d \theta ), \qquad   i: 1, \dots  m
   \\
 B_{XX}(t) hk  &=  \Big(  b^j_{x x}(t) h_0 k_0 +   \int_{-\mathtt{d}}^0 b^j_{xy}(t)f_b(\theta)  h_1(\theta) k_0\, d \theta  +   \int_{-\mathtt{d}}^0 b^j_{yx}(t)f_b(\theta)  k_1(\theta)  h_0 \, d\theta  \\ & + \int_{-\mathtt{d}}^0 \int_{-\mathtt{d}}^0 b^j_{yy}(t)f_b(\theta _1) f_b(\theta _2)   h_1(\theta_1)
 k_1(\theta_2)\, d \theta_1 \, d \theta_2 \Big ) _{1 \leq j \leq d}  \\
 \Sigma_{XX}^i(t) hk  &=   \Big( \sigma^{j,i}_{x x}(t) h_0 k_0 +   \int_{-\mathtt{d}}^0 \sigma^{j,i}_{xy}(t)f_b(\theta)  h_1(\theta) k_0\, d \theta  +   \int_{-\mathtt{d}}^0 \sigma^{j,i}_{yx}(t)f_b(\theta)  k_1(\theta)  h_0 \, d\theta  \\ & + \int_{-\mathtt{d}}^0 \int_{-\mathtt{d}}^0 \sigma^{j,i}_{yy}(t)f_b(\theta _1) f_b(\theta _2)   h_1(\theta_1)
 k_1(\theta_2)\, d \theta_1 \, d \theta_2 \Big)_{1 \leq j \leq d },  \qquad   i: 1, \dots  m.\end{align*}
 Moreover we set  
and for every $h = (h_0, h_1(\cdot)), k =  (k_0, k_1(\cdot)) \in H$\begin{align*}
 H_{X}(T) h  &= h_x(T) h_0 +  \int_{-\mathtt{d}}^0 h_y(t)f_h(\theta)  h_1(\theta)\, d \theta; \quad
  L_X(t) h = \ell_x(t) h_0 +  \int_{-\mathtt{d}}^0 \ell_y(t) f_\ell(\theta)  h_1(\theta)\, d \theta  \\
 H_{XX}(T) hk  &=  h_{x x}(T) h_0 k_0 +   \int_{-\mathtt{d}}^0 h_{xy}(T)f_h(\theta)  h_1(\theta) k_0\, d \theta  +   \int_{-\mathtt{d}}^0 h_{yx}(T)f_h(\theta)  k_1(\theta)  h_0 \, d\theta  \\ & + \int_{-\mathtt{d}}^0 \int_{-\mathtt{d}}^0 h_{yy}(T)f_h(\theta _1) f_h(\theta _2)   h_1(\theta_1)
 k_1(\theta_2)\, d \theta_1 \, d \theta_2; \\
 L_{XX}(t) hk  &=  \ell_{x x}(t) h_0 k_0 +   \int_{-\mathtt{d}}^0 \ell_{xy}(t)f_\ell(\theta)  h_1(\theta) k_0\, d \theta  +   \int_{-\mathtt{d}}^0 \ell_{yx}(t)f_\ell(\theta)  k_1(\theta)  h_0 \, d\theta  \\ & + \int_{-\mathtt{d}}^0 \int_{-\mathtt{d}}^0 \ell_{yy}(t)f_\ell(\theta _1) f_\ell(\theta _2)   h_1(\theta_1)
 k_1(\theta_2)\, d \theta_1 \, d \theta_2.
  \end{align*}
Notice that, in view of Hypotheses \ref{ipotesi}, \ref{ipotesi-costo}, \ref{ip-densita-eqstato} and \ref{ip-densita-costo}, we have that $B_X,\Sigma^i_X \in L^\infty _\Fcal (\Omega \times [0,T]; L(H; \R^d))$, $B_{XX},\Sigma^i_{XX} \in L^\infty _\Fcal (\Omega \times [0,T]; L(H \times H; \R^d))$,
$L_X  \in  L^\infty _\Fcal (\Omega \times [0,T]; L(H; \R))$, $L_{XX}  \in  L^\infty _\Fcal (\Omega \times [0,T]; L(H \times H; \R))$, $H_{X} \in L^\infty(\Omega, \Fcal_T; L(H;\R))$ and $H_{XX} \in L^\infty(\Omega, \Fcal_T; L(H \times H; \R))$.

We have that the following holds:
\begin{theorem}\label{teo:var-stato}
Under Hypotheses \ref{ipotesi} and \ref{ip-densita-eqstato} equations (\ref{equazione_varI}) and (\ref{equazione_varII}) admit unique solutions, given 
respectively by $ Y^\e \in   L^2_{\mathcal{F}}(\Omega;C([0,T]; H))$ and $ Z^\e \in   L^2_{\mathcal{F}}(\Omega \times [0,T]; H)$.
Moreover, we have
\begin{equation}\label{eq:approx-state}
\E \sup_{t \in [0,T]} |X(t)-X^\e (t)|^{2} = O(\e), 
\end{equation}
\begin{equation}\label{eq:stima_varI}
\E \sup_{t \in [0,T]} |Y^\e (t)|^{2} = O(\e), 
\qquad
\E \sup_{t \in [0,T]} |Z^\e (t)|^{2} = O(\e), 
\end{equation}\vspace{-0.2truecm}
\begin{equation}\label{eq:approx-stateI}
\E \sup_{t \in [0,T]} |X^\e(t)-X(t)-Y^\e (t)|^{2} = O(\e^{2}), 
\end{equation}\vspace{-0.2truecm}
\begin{equation}\label{eq:approx-stateII}
\E \sup_{t \in [0,T]} |X^\e(t)-X(t)-Y^\e (t)-Z^\e(t)|^{2} = o(\e^{2}), 
\end{equation} 
Moreover if also Hypotheses \ref{ipotesi-costo} and \ref{ip-densita-costo} hold true, we get following expansion of the cost functional\vspace{-0.2truecm}
\begin{align}\label{costo-espansione}\vspace{-0.2truecm}
J(u^\e) &=J(u)+\E\, H_X(T)(Y^\e(T)+Z^\e(T))\>+\frac{1}{2}\E  \, 
 H_{XX} (T)(Y^\e(T))^2\\ \nonumber&+
\E\, \int_0^T\left[ \langle  L_X(t), Y ^ \e(t) + Z ^ \e(t) \rangle \right] \, dt 
+\frac{1}{2}\E \int_0^T L_{XX}(t)(Y^\e(t))^2\, dt 
+\E\int_0^T
\delta L (t)\,dt+o(\e)
\end{align}
\end{theorem}
\dim Notice that due to the abstract reformulation, we arrive to equation \eqref{eq-astr}, which is an infinite dimensional evolution equation without delay, and with finite dimensional noise. So
all the statements follows essentially from \cite [Proposition 4.3, Proposition 4.4] {HuFuTessMax}, where the infinite dimensional case is considered. 
\qed
\subsection{First order adjoint equation}\label{subsec-first}
In this subsection we will introduce the first order adjoint equation involved in the stochastic maximum principle. Such an adjoint equation turns out to be an infinite dimensional backward stochastic differential equation, see \cite{HuFuTessMax}. 
The first order adjoint equation is given by the following BSDE with values in $H$:\vspace{-0.2truecm}
\begin{align}\label{eq:adjoint-I}
  -dp(t)& = \left[ A^*p(t)+ B^*_{X} (t) G^* p(t) -  L_{X}(t)  + \sum_{i=1}^m \Sigma^{i}_X(t)^* G^* q^i(t) \right]\,dt
 +q(t)\,dW(t) \\ \nonumber 
   p(T) &= -H_X(T).
  \end{align}
  where $q(t) = (q^i(t))_{1 \leq i \leq m}$, with $q^i(t) \in H$.  
Then, by \cite{HuPeng} we have:
\begin{theorem}\label{esistenzaI}
The equation  \eqref{eq:adjoint-I} admits a unique mild solution $(p,q)$ that is a couple of processes such that $ (p,q) \in L^2_{\mathcal{F}}(\Omega; C([0,T];H)) \times L^2_{\mathcal{F}}(\Omega\times [0,T];L (\R^m;H))$ such that for  $t \in [0,T]$ and $  \mathbb{P}$- a.s.\vspace{-0.2truecm}
\begin{align}\label{Iadjoint-mild}
    p(t) = & - e^{(T-t)A^*}H_X(T) +  \int_t^T e^{(s-t)A ^*} B_X^*(s) G^* p(s) \, ds -  \int_t^T e^{(s-t)A ^*}  L_X(s) \, ds  \\ \nonumber
    & + \sum_{i=1}^m  \int_t^T e^{(s-t)A ^*}\Sigma^i_X(s)^* G^* q^{i}(s)  \, ds  + \sum_{i=1}^m\int_t^T  e^{(s-t)A ^*}  q^{i}(s)  \, d W^i(s) .
\end{align}
\end{theorem}
\dim
See \cite[Theorem 3.1]{HuPeng}

\qed

Making use of this first adjoint equation, it turns out that the expansion of the cost in formula (\ref{costo-espansione}) can be rewritten in terms of $(p,q)$.

\begin{proposition}\label{espansioneI}
Under Hypotheses \ref{ipotesi}, \ref{ipotesi-costo}, \ref{ip-densita-eqstato} and \ref{ip-densita-costo}, , the so called \textit{duality relation} holds true
 \begin{align}\label{duality relation}
\E p(T)Y^\e(T)&=\E p_0(T)Y_0^\e(T) +  \E p_1(T)Y^\e_1(T)  = \\ \nonumber& =\E\int_0^T \left(L_X(t)Y^\e(t)+\sum_{i=1}^m  [\delta \Sigma ^{i} (t) ^* G^*q^i(t)] \right) I_{E_\e}(t) \,dt
\end{align}
from which the following expansion for the cost holds:\vspace{-0.2truecm}
 \begin{align}\label{costo-espansione-I}
& J(u)-J(u^\e) = 
\E\int_0^T \left(-
\delta L (t)+ \delta B^*(t) G^* p(t)
 +  \sum_{i=1}^m  [\delta \Sigma ^{i} (t) ^* G^*q^i(t)] \right) I_{E_\e}(t)\, dt
\\ &\nonumber
-\frac{1}{2} \E \,  H_{XX}(X(T))Y^\e(T)^2 
-\frac{1}{2}\E \int_0^T L_{XX}(t)Y^\e(t)^2\, dt \\ \nonumber&
+\frac{1}{2}\E \int_0^T B^*_{XX}(t) G^* p(t)Y^\e(t)^2  \, dt 
+\frac{1}{2}\E \int_0^T \sum_{i=1}^m   [\Sigma ^{i}_{XX}(t)^*  G^* q^i(t)Y^\e(t)^2 ]\, dt 
 +o(\e)
\end{align}
\end{proposition}
\dim 
This expansion  follows from the computation of $d\<p(t), Y^\e(t)\>$ and $ d\<p(t), Z^\e(t)\>$ by applying the It\^o formula. Notice that the differentials $dp(t),\,dY^\e(t),\,dZ^\e(t)$ are not well defined since $A$ and $A^*$ are unbounded operators. In order to work with the differentials we have to approximate $A$ and $A^*$ with their Yosida approximants $A_\lambda$ and $A^*_\lambda$.
This kind of procedure is standard when dealing with calculus in infinite dimensional Hilbert spaces, for details we refere, for instance to \cite[Theorem 6]{GMO} where a similar problem is treated with more details.
\qed

\subsection{Second order adjoint}\label{sec-second}
We introduce now the second adjoint equation that is a BDSE with values in $S^2(H)$, the Hilbert space of Hilbert-Schmidt symmetric operators from $H$ to $H$, see \cite{GuaTess}, \cite{FuHuTe-CR}, \cite{GuaTess1}.
\begin{equation}\label{second_adjoint} \left\{
\begin{array}{rl} 
-dP (t) &= \left( P(t)A+ A^*P(t) + P(t) GB_X(t)+ B^*_X(t) G^*P(t)\right) \, dt +  \mathcal{H}_{X,X} (t) \, dt +\\ \\
& +\sum_{i
=1}^m [ 
 \Sigma^i_X(t)^ *  G^* P (t) G \Sigma_X^i(t)]\, dt + \sum_{i=1}^m \left( \Sigma^i_X(t) ^*  G^* Q^i (t) +  
 Q^i(t)G \Sigma^i_X(t)\right)\,d t   \\ \\ 
& + \sum_{i=1} ^m Q^{i}(t)dW^i(t)
  \vspace{1mm}\\  
   P(T) &=- H_{XX} (T).
\end{array}
\right.
\end{equation}
where\vspace{-0.2truecm}
\begin{equation}\label{hamiltonian}
  \mathcal{H}(t)= \mathcal{H}(t, p(t), q(t))  =  
 B ^*(t) G^*p(t)  +  \sum_{i=1}^m \Sigma^ {i}(t)^* G^*q^{i}(t)  - L (t) 
\end{equation}
where $(p,q) \in  L^2_{\mathcal{F}}(\Omega; C([0,T];H)) \times L^2_{\mathcal{F}}(\Omega\times [0,T];L (\R^m;H)) $ is the 
 mild solution to \eqref{eq:adjoint-I}. Notice that the Hamiltonian function $\calh$ depends also on the state $X$ and on the control $u$:
 \begin{equation*}
     \mathcal{H}(t)=\mathcal{H}(t,X(t),u(t),  p(t), q(t)):
 \end{equation*}
for brevity we avoid writing this dependence unless it is necessary.
We notice that:
\begin{lemma} Under Hypotheses
\ref{ipotesi-costo} and \ref{ip-densita-costo}, we have that $H_{XX} \in L^\infty_{\mathcal{F}}(\Omega\times [0,T];S^2(H))$ and  $\mathcal{H}_{XX} \in L^2_{\mathcal{F}}(\Omega\times [0,T]; S^2(H))$.
 \label{regolarita dati} 
\end{lemma}\vspace{-0.2truecm}
\dim
We will show that $\Sigma ^i _{XX}(\cdot)^* G^* q^i \in  L^2_{\mathcal{F}}(\Omega\times [0,T]; S^2(H))$,  the other terms  will follow in the same way. 
We have that:\vspace{-0.5truecm}
\begin{align*}
   \Sigma ^i _{XX}(t)^* G^* q^i (t)= \sum_{j=1}^d  \Sigma ^{j,i} _{XX} (t) q^{i,j}_0 (t)\vspace{-0.6truecm}
\end{align*}
where $q^i (t) \in H $ and  $q^{i} (t)= (q^{i,j}_0 (t),q^{i,j}_1 (t))  \in  \R^d \oplus L^2((-d,0); \R^d)$ for every $i=1, \dots , m$. We notice that $G^*q^i(t)=(q^{i,j}_0(t))_{j}$, since $G^*$ is the projection on the first component of $H$.
Now, let us consider a complete orthonormal basis $\{ e^k\}_{k \geq 1}$ of $H= \R^d \oplus L^2((-d,0); \R^d)$, where $e^k=(e^k_0,0), \,k=1,...,d$ with $e_0^k, \,k=1,...,d$ a complete orthonormal basis in $\R^d$ and $e^k=(0,e^{k-d}_1), \,k\geq d+1$, where $e^{j}_1, \,j\geq 1$ is a complete orthonormal basis in $L^2((-d,0); \R^d)$.
 Hence we get that, for any $i,j$:\vspace{-0.2truecm}
 \begin{align*}
  \sum_{k \geq 1} | \Sigma ^{j,i} _{XX} (t)q^{i,j}_0(t)   e_k| ^2 & =   \sum_{k= 1} ^d |e^k_0|^2  |q^{i,j}_0 (t)|^2  |\sigma ^{j,i}_{xx} (t)|^2_{\R^{d \times d}} +  \sum_{k \geq d+1}  |\sigma ^{j,i}_{xy} (t)|^2_{\R^{d \times d}}  |q^{i,j}_0 (t)|^2
  \Big| \int_{-\mathtt{d}}^0 f_\sigma (\theta) e_1^k(\theta) \, d\theta\Big|^2 \\
  & +  \sum_{k= 1} ^d |\sigma ^{j,i}_{xy} (t)|^2_{\R^{d \times d}}  |q^{i,j}_0 (t)|^2
  |e^k_0|^2 \int_{-\mathtt{d}}^0 f^2_\sigma (\theta) \, d\theta \\ & +  \sum_{k \geq d+1}  \int_{-\mathtt{d}}^0 f ^2_\sigma(\theta_2) |\sigma ^{j,i}_{yy} (t)|^2_{\R^{d \times d}}  |q^{i,j}_0 (t)|^2
  \Big| \int_{-\mathtt{d}}^0 f_\sigma (\theta_1) e_1^k(\theta_1) \, d\theta_1\Big|^2 \, d \theta_2 \\
  &\leq C(n) |q^{i,j}_0 (t)|^2  (1 +  ||f||_{L^2}^4 )\vspace{-0.2truecm}
 \end{align*}
\qed\vspace{-0.2truecm}
\begin{theorem}\label{teo-exist-IIorder}
Under Hypotheses \ref{ipotesi}, \ref{ipotesi-costo}, \ref{ip-densita-eqstato} and \ref{ip-densita-costo}, the second order adjoint equation \eqref{second_adjoint} admits a unique mild solution  that is two couple of processes $(P,Q)$ such that 
$ P\in L^2_{\mathcal{F}}(\Omega; C([0,T]; S^2(H))) $  and $
Q \in L^2_{\mathcal{F}}(\Omega\times [0,T]; L(\R^m;S^2(H)))$.
\end{theorem}
\dim By Lemma \ref{regolarita dati} we have enough regularity to apply \cite [Theorem 5.4]{GuaTess}, thus equation \eqref{second_adjoint} has a unique mild solution $(P,Q)$ such that
$P\in L^2_{\mathcal{F}}(\Omega; C([0,T];S^2(H))$ and 
$Q \in L^2_{\mathcal{F}}(\Omega\times [0,T]; L(\R^m;S^2(H)))$. 
\qed

\smallskip
Thus we have the following reduction for the cost expansion, that gets rid of the quadratic terms in \eqref{costo-espansione-I}

\begin{proposition}
Under Hypotheses \ref{ipotesi}, \ref{ipotesi-costo}, \ref{ip-densita-eqstato} and \ref{ip-densita-costo},  the following expansion for the cost holds true:\vspace{-0.2truecm}
 \begin{align}\label{costo-espansione-II}
 J(u)-J(u^\e) = &
\E\int_0^T \left(-
\delta L (t)+ \delta B^*(t) G^* p(t)
 + \frac{1}{2}  \sum_{i=1}^m [\delta \Sigma ^{i} (t) ^* G^*q^i(t)] \right. \\ \nonumber&
+\frac{1}{2}  \sum_{i=1}^m {\tr}[  \Sigma^i(t)^*G^*P(t)G\delta \Sigma^i(t)] I_{E_\e}(t)\, dt
 +o(\e)
\end{align}
\end{proposition}
\dim
From \eqref{costo-espansione-I} in Proposition \ref{espansioneI} and  taking into account the definitions of $\calh$ we get:\vspace{-0.2truecm}
 \begin{align}\label{costo-espansione-I-calh}
& J(u)-J(u^\e) = 
\E\int_0^T \left(-
\delta L (t)+ \delta B^*(t) G^* p(t)
 +  \sum_{i=1}^m  [\delta \Sigma ^{i} (t) ^* G^*q^i(t)] \right) I_{E_\e}(t)\, dt
\\ &\nonumber
-\frac{1}{2} \E \,  H_{XX}(X(T))Y^\e(T)^2 
+\frac{1}{2}\E \int_0^T \calh_{XX}(t)Y^\e(t)^2\, dt  +o(\e)
\end{align}
Here and in the following it will be useful to write $P(t)$ in a matrix form, namely
\begin{equation}\label{P-forma-matriociale}
P(t)=\left(\begin{array}{cc}
     P_{0,0}(t)&P_{0,1} (t) \\
 P_{1,0}(t)&P_{1,1}(t)
\end{array}
\right) \text{ where } P_{0,1} (t) =P_{1,0}(t)^*,
\end{equation}
and $P_{0,0}(t)\in L(\R^d,\R^d)\simeq \R^{d\times d}$, $P_{1,0}(t)\in L(\R^d,L^2([-\mathtt{d},0],\R^d))\simeq L^2([-\mathtt{d},0],\R^d)\times ....\times L^2([-\mathtt{d},0],\R^d)$, $P_{1,1}(t)\in S^2(L^2([-\mathtt{d},0],\R^d),L^2([-\mathtt{d},0],\R^d)$ and are all symmetric operators. 
Hence we compute 
\begin{align}\label{traccia}
  P(t) Y^\e(t)^2 &=\<P(t)Y^\e(t), Y^\e(t)\>:=\sum_{i,j=1}^d P^{i,j}_{0,0}(t)Y^i_{0}(t)Y^j_{0}(t)+
\sum_{i=1}^d \<P^{i}_{1,0}(t)Y^i_{0}(t),Y_{1}(t)\>_{L^2([-\mathtt{d},0],\R^d)}\\ \nonumber
&+\sum_{i=1}^d \<Y_{1}(t),P^{i}_{1,0}(t)Y^i_{0}(t)\>_{L^2([-\mathtt{d},0],\R^d)}+ \<P_{1,1}(t)Y_{1}(t),Y_{1}(t)\>_{L^2([-\mathtt{d},0],\R^d)}:
\end{align}
by applying the It\^o formula we have \vspace{-0.2truecm}
\begin{align*}
     d\<P(t)Y^\e(t), Y^\e(t)\> &= \sum_{i=1}^m\<P(t)G\delta \Sigma^i(t),G\delta \Sigma^i(t)\>-\<\calh_{XX}(t)Y^\e(t), Y^\e(t)\>\\&=
\sum_{i=1}^m \Sigma^i(t)^*G^*P(t)G\delta \Sigma^i(t)  -\calh_{XX}(t)Y^\e(t)^2 
\end{align*}
so integrating between $0$ and $T$ and taking the expectation we get 
\begin{align}\label{relaggiunta II}
  &- \E H_{XX}(T)Y^\e(T)^2 + \E\int_0^T \calh_{XX}(t)Y^\e(t)^2 \, dt = 
      \E\int_0^T \sum_{i=1}^m{\tr}[  \Sigma^i(t)^*G^*P(t)G\delta \Sigma^i(t)]\,dt
\end{align}
Substituting into \eqref{costo-espansione-I-calh} we obtain the desired expansion for the cost.\qed

\subsection{Maximum principle}\label{max-princ}


Now we are able to prove a version of the Stochastic Maximum Principle, in its necessary form, for the control problem with state equation and cost functional given
by \eqref{equazione_stato} and \eqref{costo}, respectively.
\noindent We recall that the expansion of the cost is given in (\ref{costo-espansione}) and it is rewritten in (\ref{costo-espansione-II}) in terms of the pair of processes
$(p,q)$ solution of the first adjoint equation (\ref{eq:adjoint-I}) and of $(P,Q)$ solution of the second order adjoint equation \eqref{second_adjoint}.

In the following Theorem the stochastic maximum principle is written in terms of the Hamiltonian function $\calh$ defined in (\ref{hamiltonian}).
\begin{theorem}\label{teoMaxPrinc}
Let \ref{ipotesi}, \ref{ipotesi-costo}, \ref{ip-densita-eqstato} and \ref{ip-densita-costo}  be satisfied and suppose that $(X,u)$ is an optimal pair for the control problem, 
and let us consider $u^\e$ defined in (\ref{spike_def}).
There exist two pairs of processes $ (p,q) \in L^2_{\mathcal{F}}(\Omega; C([0,T];H)) \times
 L^2_{\mathcal{F}}(\Omega\times [0,T];L(\R^m,H))$, solution to the first order adjoint equation \eqref{eq:adjoint-I}, 
and $(P,Q)\in L^2_{\mathcal{F}}(\Omega; C([0,T];S^{2}(H))) \times \in L^2_{\mathcal{F}}(\Omega\times [0,T];L(\R^{m};S^{2}(H))$, solution of the second order adjoint equation \eqref{second_adjoint}, such that the following variational inequality holds 
 \begin{align*}\label{eq:var-ham}
  &\calh(t, X(t),  u(t),p(t),q(t))-   \calh(t, X(t), v,p(t),q(t)) \\ \vspace{-0.2truecm} &+\frac{1}{2} \sum_{i=1}^m (\delta \Sigma^i(t))^*  G^* P(t) G \delta  \Sigma^i(t)\leq0
  \quad \forall\,v\in U,\quad \text{a.e. } t\in [0,T],\, \P-\text{a.s.},
 \end{align*}
 where $\calh $ is the Hamiltonian function defined in \eqref{hamiltonian}.
 Equivalently, for $ \text{ a.e. } t\in [0,T],\, \P-\text{a.s.}$,\vspace{-0.2truecm}
 \begin{align*}
  &\calh(t, X(t),  u(t),p(t),q(t)) +\frac{1}{2} \sum_{i=1}^m ( \Sigma^i(t, X(t),u(t)))^* G^* P(t) G \Sigma^i(t, X(t),u(t))\\ \vspace{-0.5truecm}
  & \;=\inf_{v\in U}  \calh(t, X(t),v,p(t),q(t))
  +\frac{1}{2} \sum_{i=1}^m (\Sigma^i(t, X(t),v)^* G^*P(t) G\Sigma^i(t, X(t),v) .\end{align*}
\end{theorem}
\dim By (\ref{costo-espansione-II}) we get that\vspace{-0.2truecm}
\begin{align*}
0&\geq J(u)-J(u^\e)
  = \E\int_0^T \left(-
\delta L (t)+ \delta B^*(t) G^* p(t)
 + \frac{1}{2} \sum_{i=1}^m  [\delta \Sigma ^{i} (t) ^* G^*q^i(t)] \right. \\ \nonumber&\qquad\qquad\left.
+\frac{1}{2}  \sum_{i=1}^m  \delta \Sigma ^{i} (t) ^* G^*P(t)G\Sigma ^{i} (t)\right) I_{E_\e}(t)\, dt
 +o(\e)
\\ \vspace{-0.2truecm}\nonumber& = 
  \E\int_0^T \big[\calh(t, X(t),  u(t),p(t),q(t))-   \calh(t, X(t), v,p(t),q(t)) +\frac{1}{2} \sum_{i=1}^m {\tr}[(\delta \Sigma^i(t))^* P(t) \delta \Sigma^i(t)]\big] I_{E_\e(t)} \, dt
   +o(\e)
\end{align*}
Notice that the right hand side is not equal to $0$ only for $t\in E_\e$. So 
if we divide both sides by $\e$ and we let $\e\rightarrow 0$ we get\vspace{-0.2truecm} 
\[
 0\geq \calh(t, X(t),  u(t),p(t),q(t))-   \calh(t, X(t), v,p(t),q(t)) +\frac{1}{2} \sum_{i=1}^m (\delta \Sigma^i(t))^*G^* P(t) G\delta \Sigma^i(t)
\]
$\forall\,v\in U$, a.e. $t\in [0,T],\, \P-$ a.s., and this concludes the proof. 
\qed
\begin{remark}\label{rm:prima-componente}
We notice that the Hamiltonian function can be written as\vspace{-0.4truecm}
\begin{align} \label{hamiltoniana ridotta}
 \mathcal H(t)&= \mathcal{H}(t, p(t), q(t))  =  
 B ^*(t) G^*p(t)  +  \sum_{i=1}^m \Sigma^ {i}(t)^* G^*q^{i}(t)  - L (t) \\  \vspace{-0.5truecm}\nonumber
&= \sum_{j=1}^d b^j(t) p^j_0(t)
 + \sum_{j=1}^d \sum_{i=1}^m\sigma^ {j,i}(t) q^{i,j}_0(t) 
-\ell(t):=\mathcal{H}^a(t)(t, p_0(t),q_0(t)).
\end{align}
Here $\mathcal H:[0,T]\times H\times L (\R^m;H)\rightarrow \R$ while $\mathcal H^a:[0,T]\times \R^d\times L (\R^m;\R^d)\rightarrow \R$, and the link between $\calh$ and $\calh^a$ is given by 
$\calh(t, p(t), q(t))=\calh^a(t, G^*p(t), G^*q(t))$.
\end{remark}

We underline the fact that in the formulation of the stochastic maximum principle only the first component $(p_0,q_0)$ of the solution of the first order equation \eqref{eq:adjoint-I} enters, and for what concerns the solution $(P,Q)$ of the second order equation \eqref{second_adjoint} only $P_{0,0}$ appears. Indeed we notice that the hamiltonian function can be written as 
in \eqref{hamiltoniana ridotta}
and the term involving $P$ can be written as 
\vspace{-0.2truecm}
\begin{equation*}
\frac{1}{2} \sum_{i=1}^m (\delta \Sigma^i(t))^*G^* P(t) G\delta \Sigma^i(t)= \frac{1}{2}\sum_{i=1}^m(\delta \Sigma^i(t))^* P_{0,0}(t)) \delta \Sigma^i(t).
\end{equation*}
This will be crucial in the next Sections when we consider in the state equation dependence on the past through general measures. We will treat such a general case using a first order adjoint equation which is an anticipated BSDE, then we rewrite the equations for the first and second variation of the state as abstract evolution equations in $H$, we approximate the general measures by means of absolutely continuous measures and we write an approximating second order adjoint equation.
When we will let $n\rightarrow\infty$, we will study only the convergence of the process appearing in the formulation of the stochastic maximum principle.

\section{The case of a general measure}\label{sec-gen-meas}

We consider now a more general case where the measures involved are finite regular measures on $[-\mathtt{d},0]$ ( see e.g. \cite{Dud} for the definition of regular measure ) not necessarily absolutely continuous with respect to $\lambda$, the Lebesgue measure on $[-\mathtt{d},0]$. Notice that in this case the Hilbert space $H=\mathbb{R}^n\oplus L^{2}\left(\left[-\mathtt{d},0\right] ,\mathbb{R}^n\right)$ where we reformulate the problem should be replaced by the Banach space $E=\mathbb{R}^n\oplus C\left(\left[-\mathtt{d},0\right] ,\mathbb{R}^n\right)$, and the dual processes $(p,q)$ would take values in $E^*$, moreover the second order adjoint equation would take values in a subspace of $L(E^*, E^*)$.

We do not reformulate directly the problem as an infinite dimensional problem without delay, but we consider first and second variation of the state and their (first order) adjoint equation which turns out to be an anticipated BSDE, see \cite{ElYa}, \cite{PengYang}, and \cite{GuaMas} where the case of dependence on the past trajectory through a general measure is treated.
\newline We underline the fact that from now on we assume that the final cost does not depend on the past, that is we assume that $h: \R^d\rightarrow R$. Indeed, in order to consider a general dependence on the past we should work with a new form of anticipated BSDEs, see \cite{GuaMas}, and this adds new difficulties.

\subsection{First and second order variation of the state} \label{sec:var.state.gen}

We consider the stochastic delay differential equation \eqref{equazione_stato}, without reformulating it as an abstract evolution equation in $H$ and we write the first and second variation. So the equations for the first and second variation are themselves stochastic (linear) delay differential equations.
In the following we use the notation\vspace{-0.2truecm}
\begin{align*}
&\psi(t):=\psi(t,x(t),\dis\int_{-\mathtt{d}}^0 x(t+\theta )\mu_\psi(d\theta),u(t)), \quad \psi_\rho(t)= \psi_\rho(t,\dis\int_{-\mathtt{d}}^0 x(t+\theta )\mu_\psi(d\theta),u(t))\\
&\psi^\e(t):= \psi(t,\dis\int_{-\mathtt{d}}^0 x(t+\theta )\mu_\psi(d\theta),u^\e(t)), \\
&\delta \psi(t):= \psi(t,\dis\int_{-\mathtt{d}}^0 x(t+\theta )\mu_\psi(d\theta),u^\e(t))- \psi(t,\dis\int_{-\mathtt{d}}^0 x(t+\theta )\mu_\psi(d\theta),u(t)),
\end{align*}
where $\psi= b ,\,\sigma ,\, \ell$ and $\rho=x,y,u$.
\newline The equation for the first variation of the state is the following, see also \cite{ChenWu},\vspace{-0.2truecm}
\begin{equation}\label{agg.prima.fin.dim}
\left\{ \begin{array}{ll}
dy^\e(t) =  [\displaystyle b_x(t) y^\e(t) + b_y(t) \int_{-\mathtt{d}}^ 0y^ \e (t+ \theta) {\mu} _b(d\theta) ] \, dt + [\sigma_x(t) y^\e(t) \\ \\  \qquad\quad \,  + \sigma_y(t) \displaystyle\int_{-\mathtt{d}}^0 y^ \e (t+ \theta) \mu_{\sigma}(d\theta) ]\,dW(t)
   +   \, \delta {\sigma}(t) ]\, d W(t), \quad    t \in [0,T], \\ \\
  y^\e(0)=0, \quad y^\e(\theta)=0, \quad-\mathtt{d} \leq \theta <0.
\end{array} \right.
\end{equation}
and the second variation is, for $ t \in [0,T]$,
\begin{equation}\label{agg.seconda.fin.dim}
\left\{ \begin{array}{ll}
dz^\e(t) =  \displaystyle \Big[b_x(t) z^\e(t) +\ b_y(t) \int_{-\mathtt{d}}^ 0z^ \e (t+ \theta) {\mu} _b(d\theta) + \, \delta {b}(t) \,\big] dt +  \displaystyle \frac{1}{2}\Big[b_{xx}(t)y^\e(t)y^\e(t)\\ \qquad\quad \,  + \displaystyle\int_{-\mathtt{d}}^ 0 \displaystyle\int_{-\mathtt{d}}^ 0b_{yy}(t) y^\e(t+\theta_1)y^\e(t+\theta_2){\mu} _b(d\theta_1) {\mu} _b(d\theta_2)
+2 \displaystyle\int_{-\mathtt{d}}^0 b_{xy}(t)y^\e(t)  y^{\e}(t+\theta)\mu_b(d\theta)\Big] \, dt\\ 
\qquad\quad \,  \displaystyle   +\Big[\sigma_x(t) z^\e(t) +{ \sigma}_y(t) \int_{-\mathtt{d}}^ 0z^ \e (t+ \theta) {\mu}_b (d\theta)+\delta {\sigma}_x(t) y^\e(t) +\delta {\sigma}_y(t) \displaystyle \int_{-\mathtt{d}}^ 0y^ {\e} (t+ \theta) \mu_\sigma(d\theta)\Big] \,d W(t)\\ 
\qquad\quad \,  \displaystyle + \frac{1}{2}\Big[\sigma_{xx}(t)y^\e(t)y^\e(t) +\int_{-\mathtt{d}}^ 0\int_{-\mathtt{d}}^ 0\sigma_{yy}(t) y^\e(t+\theta_1)y^\e(t+\theta_2) {\mu} _\sigma(d\theta_1) {\mu} _\sigma(d\theta_2)\Big] 
  \, d W(t)\\
 \qquad\quad \,  +\Big[\displaystyle\int_{-\mathtt{d}}^0\sigma_{xy}(t)y^\e(t)  y^\e(t+\theta){\mu}_\sigma(d\theta)\Big] 
  \, d W(t) \\ 
  z^\e(0)=0, \quad z^\e(\theta)=0, \quad-\mathtt{d} \leq \theta <0.
\end{array} \right.
\end{equation}
We recall that, see e.g. \cite{Moha}, under Hypothesis \ref{ipotesi}, we have that equations \eqref{agg.prima.fin.dim}, \eqref{agg.seconda.fin.dim} admit a unique solution such that $\forall\,p\geq 1$:\vspace{-0.2truecm}
  \begin{align}\label{stima sup y}
        & \E \sup_{t\in [-\mathtt{d},T]} |y^{\e}(t)|^p < +\infty,\qquad\E \sup_{t\in [-\mathtt{d},T]} |z^{\e}(t)| ^p < +\infty, \vspace{-0.2truecm}
      \end{align} 
Moreover \vspace{-0.3truecm}
\begin{equation}\label{eq:stima_varI-concreto}
\E \sup_{t \in [0,T]} |y^\e (t)|^{2} = O(\e), 
\qquad
\E \sup_{t \in [0,T]} |z^\e (t)|^{2} = O(\e), 
\end{equation}\vspace{-0.5truecm}
\begin{equation}\label{eq:approx-stateI-concreto}
\E \sup_{t \in [0,T]} |x^\e(t)-x(t)-y^\e (t)|^{2} = O(\e^{2}), 
\end{equation}\vspace{-0.5truecm}
\begin{equation}\label{eq:approx-stateII-concreto}
\E \sup_{t \in [0,T]} |x^\e(t)-x(t)-y^\e (t)-z^\e(t)|^{2} = o(\e^{2}), 
\end{equation} 

\subsection{Duality relation I: first order adjoint equation}\label{subsec-dualeI-anticp}
 We consider the following anticipated backward stochastic differential equation (ABSDE) (see \cite{PengYang} and \cite{YangEll} for more details on ABSDEs)
\begin{equation}\label{anticipante-sec}
 \left\lbrace\begin{array}{l}
{p}^a(t)=\dis\int_t^T\ell_x(s)ds + \int_t^T \int_{-\mathtt{d}}^0 \ell_y(s-\theta) \mu_{\ell}(d\theta)\,ds \\
\dis\int_t^T\dis b_x(s) {p}^a(s) ds  +
\dis\int_t^T\dis \int_{-\mathtt{d}}^0 \, \E^{\Fcal_s}b_y(s-\theta)p^a(s-\theta)\mu_b(d\theta)\, ds \,   \\ 
+\dis\int_t^T \dis \sigma_x(s) {q}^a(s) ds+\dis\int_t^T \dis \int_{-\mathtt{d}}^0 \E^{\Fcal_s}\sigma_y(s-\theta) {q}^a(s-\theta)\mu_\sigma(d\theta) \,  ds  -
\dis\int_t^T  q^a(s)dW(s) \displaystyle + h_x(T),  \\ \\ 
 {p}^a(T-\theta), \; {q}^a(T-\theta)=0 \quad \forall \,\theta \in [-\mathtt{d},0).
 \end{array}
 \right.
 \end{equation}
 where $(p^a,q^a) \in L^2_{\mathcal{F}}(\Omega; C([0,T+\mathtt{d}];\R^d)) \times L^2_{\mathcal{F}}(\Omega\times [0,T+\mathtt{d}];L (\R^m;\R^d))$. Equation \eqref{anticipante-sec} turns out to be the first order adjoint equation of \eqref{agg.prima.fin.dim}. 
\begin{proposition}\label{prop:Ito}
Under  Hypotheses \ref{ipotesi} and \ref{ipotesi-costo} the  following  relations hold: \vspace{-0.2truecm}
       \begin{align}\label{duality relation y}
\E p^a(T)y^{\e}(T))=\E h_x(T)y^{\e}(T)  &=\E\dis\int_0^T[\ell_x(s),y^\e(s)\>+q^a(s)\delta\sigma(s) I_{E_\e}(s)]\,ds
\\\nonumber&+
\E\dis\int_0^T\int_{-\mathtt{d}}^0 \ell_y(s)   y^\e(s+\theta)  \mu_\ell(d\theta)  \, ds 
\end{align}\vspace{-0.6truecm}
     \begin{align}\label{duality relation z}
&\E p^a(T)z^{\e}(T))=\E h_x(x(T)) z^\e(T)=\E\dis\int_0^T[\ell_x(s),z^\e(s)+ \int_{-\mathtt{d}}^0 \ell_y(s)   z^\e(s+\theta)  \mu_\ell(d\theta) ] \, ds  \\ \nonumber&+ \E\dis\int_0^Tq^a(s)[ \delta {\sigma}_x(s) y^\e(s)+   \displaystyle  \delta {\sigma}_y(s) \displaystyle \int_{-\mathtt{d}} ^ 0 y^ {\e} (s+ \theta) \mu_\sigma(d\theta)) I_{E_\e}(s)]\,ds\\ \nonumber
&\quad+ \E\dis\int_0^T\frac{1}{2}p^a(s)\left[ b_{xx}(s)y^\e(s)y^\e(s)+\int_{-\mathtt{d}} ^ 0 \int_{-\mathtt{d}} ^ 0 b_{yy}(s) y^\e(s+\theta_1)y^\e(t+\theta_2){\mu} _b(d\theta_1) {\mu} _b(d\theta_2)\right]\,ds\\ \nonumber
&\quad +\E\dis\int_0^Tp^a(s)\left[\int_{-\mathtt{d}}^0 b_{xy}(s)y^\e(s)  y^{\e}(s+\theta)\mu_b(d\theta)+\delta b(s) I_{E_\e}(s)\right]\,ds\\ \nonumber
 &\quad+ \E\dis\int_0^T\frac{1}{2}q^a(s)\left[\sigma_{xx}(s)y^\e(s)y^\e(s)+\int_{-\mathtt{d}}^ 0 \int_{-\mathtt{d}}^ 0 \sigma_{yy}(s) y^\e(s+\theta_1)y^\e(t+\theta_2){\mu} _\sigma(d\theta_1) {\mu} _\sigma(d\theta_2)\right]\,ds\\ \nonumber
&\quad +\E\dis\int_0^Tp^a(s)\int_{-\mathtt{d}}^0 \sigma_{xy}(s)y^\e(s)  y^{\e}(s+\theta)\mu_\sigma(d\theta)\>\,ds
\end{align} 
\end{proposition}
\dim See the Appendix \ref{appendix}.

\begin{corollary}\label{espansioneI-gen}
Under  Hypotheses \ref{ipotesi} and \ref{ipotesi-costo}, the following expansion for the cost holds true:
 \begin{align}\label{costo-espansione-I-gen}
& J(u)-J(u^\e) = 
\E\int_0^T \left(-
\delta \ell (t)+ p^a(t)\delta b(t) +q^a(t)\delta \sigma(t) \right) I_{E_\e}(t)\,dt \\ 
 &\nonumber-\frac{1}{2} \E \<h_{xx}(y(T))y^\e(T),y^\e(T)\> 
-\frac{1}{2}\E \int_0^T \<\ell_{xx}(t)y^\e(t),y^\e(t)\>\, dt 
\\ \nonumber& 
-\frac{1}{2}\E \int_0^T \int_{\mathtt{d}}^ 0\int_{\mathtt{d}}^ 0\ell_{yy}(t) y^\e(t+\theta_1)y^\e(t+\theta_2){\mu} _\ell(d\theta_1) {\mu} _\ell(d\theta_2)\,dt
-\E\dis\int_0^T\int_{-\mathtt{d}}^0 \ell_{xy}(t)y^\e(t)  y^{\e}(t+\theta)\mu_\ell(d\theta)\,dt
\\ \nonumber& 
+\E\dis\int_0^T\frac{1}{2}\<p^a(t), b_{xx}(t)y^\e(t)y^\e(t)+\int_{-\mathtt{d}}^ 0\int_{-\mathtt{d}}^ 0b_{yy}(t) y^\e(t+\theta_1)y^\e(t+\theta_2){\mu} _b(d\theta_1) {\mu} _b(d\theta_2)\>\,ds\\  \nonumber
&\quad +\E\dis\int_0^T\<p^a(t),\int_{-\mathtt{d}}^0 b_{xy}(t)y^\e(t)  y^{\e}(t+\theta)\mu_b(d\theta)\>\,dt\\  \nonumber
 &\quad+ \E\dis\int_0^T\frac{1}{2}\<q^a(t), \sigma_{xx}(t)y^\e(t)y^\e(t)+\int_{-\mathtt{d}}^ 0\int_{-\mathtt{d}}^ 0\sigma_{yy}(t) y^\e(t+\theta_1)y^\e(t+\theta_2){\mu} _\sigma(d\theta_1) {\mu} _\sigma(d\theta_2)\>\,dt\\  \nonumber
&\quad +\E\dis\int_0^T\<q^a(t),\int_{-\mathtt{d}}^0 \sigma_{xy}(t)y^\e(t)  y^{\e}(t+\theta)\mu_\sigma(d\theta)\>\,dt
 +o(\e)
\end{align}
\end{corollary}
\dim This expansion  follows from \eqref{duality relation y} and \eqref{duality relation z} and \eqref{eq:stima_varI-concreto}
\qed

\subsection{Regular approximations of the first and second variation}

In order to get rid of the term containing square of the first order variation, that is terms containing $y^\e(t+\theta_1)y^\e(t+\theta_2),\, \theta_1,\theta_2\in [-\mathtt{d},0]$, we have to consider a second order adjoint equation. At this level we pass to an infinite dimensional reformulation, since without reformulation it would be necessary to consider many second order adjoint equations in order to get rid of the terms $y^\e(t+\theta_1)y^\e(t+\theta_2),\, \theta_1,\theta_2\in [-\mathtt{d},0]$ with $\theta_1\neq \theta_2$.

As already discussed at the beginning of Section \ref{sec-gen-meas} we are considering the case when the measures involved are finite regular measures on $[-\mathtt{d},0]$, not necessarily absolutely continuous with respect to the Lebesgue measure on $[-\mathtt{d},0]$. In order to work in the Hilbert space $H=\mathbb{R}^n\oplus L^{2}\left(\left[-d,0\right] ,\mathbb{R}^n\right)$ where we want to reformulate the problem, we proceed by an approximation procedure. First we recall the following approximation result, see \cite{GMO}.
\begin{lemma}\label{lemma:approx-measure}
 Let $\mu$ be a finite  regular measure on $[-\mathtt{d},0]$.
 There exists a sequence $(\mu^n)_{n\geq 1}$ of finite regular measure on $[-\mathtt{d},0]$, absolutely continuous with respect to $\lambda$, such that\vspace{-0.2truecm}
 \begin{equation}\label{conv-a-mu}
  \mu=w^*-\lim_{n\rightarrow\infty}\mu^n,
 \end{equation}
that is for every $\phi\in C_b([-\mathtt{d},0];\R)$\vspace{-0.2truecm}
\begin{equation}\label{conv-a-mu:expl}
  \int_{-\mathtt{d}}^0\phi (s)\,d\mu (s)=\lim_{n\rightarrow\infty}\int_{-\mathtt{d}}^0 \phi (s)\, d\mu^n (s).
  \end{equation}
  Moreover $ |\mu|,|\mu^n| \leq  K$, for some $K>0$ independent of $n$.
\end{lemma}
\begin{remark}\label{rem dens}
Notice that we can always construct an approximating sequence $\mu^n$ whose densities are square integrable functions, see \cite[Lemma 5]{GMO} thus we will assume this property from now on.
\end{remark}
We are then interested in the control problem \eqref{equazione_stato} - \eqref{costo} where now all the measures $\mu_b, \mu_\sigma$ are  just regular measures.
In such a case we cannot follow the same procedure showed in the previous Section mainly because the terms involving derivatives with respect to the past trajectory are in general regular measures, but we cannot guarantee that these terms belong to $H$. For this reason the first and the second variations as well as the first and second adjoint equations do not evolve in $H$ and $S^2(H)$ respectively. We bypass the problem introducing a regularization method.

We will proceed with an approximation technique introducing
the following family of equations, where the regular measures $\mu_b,\,\mu_\sigma$ are replaced by two sequences of regular measures $(\mu^n_b)_n, \,(\mu^n_\sigma)_n$ with densities respectively given by $f^n_b, f^n_\sigma$:\vspace{-0.2truecm} 
\begin{equation}\label{yn}
\left\{ \begin{array}{ll}
dy^{n,\e}(t) = \Big[\displaystyle b_x(t) y^{n,\e}(t) + b_y(t) \int_{-\mathtt{d}}^ 0y^{n,\e} (t+ \theta) f^n _b (\theta)d\theta \Big] \, dt   \\ \\  \qquad\quad + \Big[\displaystyle\sigma_x(t) y^{n,\e}(t) + \sigma_y(t) \int_{-\mathtt{d}}^0 y^{n,\e} (t+ \theta) f^n_{\sigma} (\theta) \, d\theta + \, \delta {\sigma}(t)\Big ] d W(t), \quad    t \in [0,T], \\ \\
  y^{n,\e}(0)=0, \quad y^{n,\e}(\theta)=0, \quad-\mathtt{d} \leq \theta <0.
\end{array} \right.\vspace{-0.2truecm}
\end{equation}
and \vspace{-0.1truecm}
\begin{equation}\label{zn}
\left\{ \begin{array}{ll}
dz^{n,\e}(t) =  \displaystyle \Big[b_x(t) z^{n,\e}(t) +\ b_y(t) \int_{-\mathtt{d}}^ 0z^{n,\e}(t+ \theta) f^n_b(\theta)\,d\theta + \, \delta {b}(t) \,\big] dt \\
\,+  \displaystyle \frac{1}{2}\Big[b_{xx}(t)y^{n,\e}(t)y^{n,\e}(t)+\int_{-\mathtt{d}}^ 0\int_{-\mathtt{d}}^ 0b_{yy}(t) y^{n,\e}(t+\theta_1)y^{n,\e}(t+\theta_2)f^n_b(\theta_1)f^n_b(\theta_2)d\theta_1 d\theta_2\\ 
\quad+2\displaystyle\int_{-\mathtt{d}}^0 b_{xy}(t)y^{n,\e}(t) y^{n,\e}(t+\theta)f^n_b(\theta)d\theta\Big] \, dt\\ 
 \, \displaystyle   +\Big[\sigma_x(t) z^{n,\e}(t) +{ \sigma}_y(t) \int_{-\mathtt{d}}^ 0z^{n,\e} (t+ \theta) f^n_b (\theta)d\theta+\delta {\sigma}_x(t) y^{n,\e}(t) +\delta {\sigma}_y(t) \displaystyle \int_{-\mathtt{d}}^ 0y^{n,\e} (t+ \theta) f^n_\sigma(\theta)d\theta\Big] \,d W(t)\\ \,\displaystyle + \frac{1}{2}\Big[\sigma_{xx}(t)y^{n,\e}(t)y^{n,\e}(t) +\int_{-\mathtt{d}}^ 0\int_{-\mathtt{d}}^ 0\sigma_{yy}(t) y^{n,\e}(t+\theta_1)y^{n,\e}(t+\theta_2) f^n_\sigma(\theta_1) f^n_\sigma(\theta_2)d\theta_1 d\theta_2\Big] 
  \, d W(t)\\
 \,+\Big[\displaystyle\int_{-\mathtt{d}}^0\sigma_{xy}(t)y^{n,\e}t) y^{n,\e}(t+\theta)f^n_\sigma(\theta)d\theta\Big] 
  \, d W(t) \\ 
  z^{n,\e}0)=0, \quad z^{n,\e}(\theta)=0, \quad-\mathtt{d} \leq \theta <0.
\end{array} \right.
\end{equation}
Note that in the terms $b_x,\,b_y,\,\sigma_x,\,\sigma_y$ we do not perform any approximation of the measures $\mu_b$ and $\mu_\sigma$, so equations \eqref{yn} and \eqref{zn} are not first and second variation of an approximating controlled stochastic delay differential equations. We have just built equations \eqref{yn} and \eqref{zn} as approximating equations of \eqref{agg.prima.fin.dim} and \eqref{agg.seconda.fin.dim} respectively.

\noindent We have the following  convergence result:
\begin{proposition}\label{prop.conv.yn-zn}
    Under hypotheses \ref{ipotesi}, we have that, for every $\e>0$ and every $n \geq 1$ every equations \eqref{agg.prima.fin.dim}, \eqref{agg.seconda.fin.dim}, \eqref{yn} and \eqref{zn} admits a unique solution such that $p\geq 1$:
           \begin{equation}\label{stima sup yn}
        \sup_{n\geq 1} \, \E \sup_{-\mathtt{d} \leq t \leq T} |y^{n,\e}(t)|^p < +\infty, \qquad
      \sup_{n\geq 1} \, \E \sup_{-\mathtt{d} \leq t \leq T} |z^{n,\e}(t)|^p < +\infty
    \end{equation}  
    Moreover it holds:
    \begin{equation}\label{conv y}
        \lim_{n \to \infty} \E \sup_{-\mathtt{d} \leq t \leq T} |y^{n,\e}(t)- y^{\e}(t)|^2=0, 
    \qquad
        \lim_{n \to \infty} \E \sup_{-\mathtt{d} \leq t \leq T} |z^{n,\e}(t)- z^{\e}(t)|^2=0, 
    \end{equation}
\end{proposition}
\dim
See the Appendix \ref{appendix}.

\subsection{ Reformulation of the regularized variations}
\vspace{0.2truecm}
Notice that in order to get the convergence of $y^{n,\e}$ to $y^{\e}$ and of $z^{n,\e}$ to $z^{\e}$ we have made computations on the delay equations without reformulation, since we needed  continuity properties. Now we reformulate equations 
 \eqref{yn} and \eqref{zn} as abstract evolution equations in  $H$.
 
We introduce the following operators in $H$: for each $h,k \in H$\vspace{-0.2truecm}
\begin{align*}
    \hat{B}^n(t) h  &=\left( b ^ j_x(t) h_0 +  \int_{-\mathtt{d}}^0 b ^j_y(t)f^n_b(\theta)  h_1(\theta)\, d \theta\right)_{1 \leq j \leq n} \\ 
   \hat{\Sigma}^n(t) h & = \left(  \sigma^{j,i}_x(t) h_0 +  \int_{-\mathtt{d}}^0 \sigma^{j,i}_y(t) f^n_\sigma(\theta)  h_1(\theta)\, d \theta \right)_{1 \leq j \leq d, 1 \leq i \leq m }  =  \left(  \hat{\Sigma}^{i,n} (t)h\right)_{ 1 \leq i \leq m }  \\
\delta\hat{\Sigma}^n(t)h & = \left( \delta \sigma^{j,i}_x(t) h_0 +  \int_{-\mathtt{d}}^0 \delta\sigma^{j,i}_y(t) f^n_\sigma(\theta)  h_1(\theta)\, d \theta \right)_{1 \leq j \leq d, 1 \leq i \leq m }  =  \left( \delta \hat{\Sigma}^{i,n} (t)h\right)_{ 1 \leq i \leq m }  \\
\tilde{B}^n(t) hk  &=  \Big(  b^j_{x x}(t) h_0 k_0 +   \int_{-\mathtt{d}}^0 b^j_{xy}(t)f^n_b(\theta)  h_1(\theta) k_0\, d \theta  +   \int_{-\mathtt{d}}^0 b^j_{yx}(t)f^n_b(\theta)  k_1(\theta)  h_0 \, d\theta  \\ & + \int_{-\mathtt{d}}^0 \int_{-\mathtt{d}}^0 b^j_{yy}(t)f^n_b(\theta _1) f^n_b(\theta _2)   h_1(\theta_1)
 k_1(\theta_2)\, d \theta_1 \, d \theta_2 \Big ) _{1 \leq j \leq d}  \\
 \tilde{\Sigma}^n(t) hk  &=   \Big( \sigma^{j,i}_{x x}(t) h_0 k_0 +   \int_{-\mathtt{d}}^0 \sigma^{j,i}_{xy}(t)f^n_b(\theta)  h_1(\theta) k_0\, d \theta  +   \int_{-\mathtt{d}}^0 \sigma^{j,i}_{yx}(t)f^n_b(\theta)  k_1(\theta)  h_0 \, d\theta  \\ & + \int_{-\mathtt{d}}^0 \int_{-\mathtt{d}}^0 \sigma^{j,i}_{yy}(t)f^n_b(\theta _1) f^n_b(\theta _2)   h_1(\theta_1)
 k_1(\theta_2)\, d \theta_1 \, d \theta_2 \Big)_{1 \leq j \leq d, 1 \leq i \leq m }  =  \left(  \tilde{\Sigma}^{i,n} (t)hk\right)_{ 1 \leq i \leq m }  \end{align*}
 and \vspace{-0.2truecm}
 \begin{align*}
  \tilde{L}^n(t) hk  &=  \Big(  \ell_{x x}(t) h_0 k_0 +   \int_{-\mathtt{d}}^0 \ell_{xy}(t)f^n_\ell(\theta)  h_1(\theta) k_0\, d \theta  +   \int_{-\mathtt{d}}^0 \ell_{yx}(t)f^n_\ell(\theta)  k_1(\theta)  h_0 \, d\theta  \\ & + \int_{-\mathtt{d}}^0 \int_{-\mathtt{d}}^0 \ell_{yy}(t)f^n_\ell(\theta _1) f^n_\ell(\theta _2)   h_1(\theta_1)
 k_1(\theta_2)\, d \theta_1 \, d \theta_2 \Big )   ,
 \end{align*}
 where $\mu_\ell$ has been replaced by a sequence of regular measures $(\mu^n_\ell)_n, $ with density given by $f^n_\ell$ and approximating $\mu_\ell$ in the sense of lemma \ref{lemma:approx-measure}.
 In view of hypotheses \ref{ipotesi}, \ref{ipotesi-costo}, \ref{ip-densita-eqstato} and \ref{ip-densita-costo} and remark \ref{rem dens}, we have that $\hat{B},\hat{\Sigma}^{i,n} \in L^\infty _\Fcal (\Omega \times [0,T]; L(H; \R^d))$, $\tilde{B}^n,\tilde{\Sigma}^{i,n} \in L^\infty _\Fcal (\Omega \times [0,T]; L(H \times H; \R^d))$,
 $\tilde{L}^n  \in  L^\infty _\Fcal (\Omega \times [0,T]; L(H \times H; \R))$.

We reformulate equations  \eqref{yn} and \eqref{zn}, see \cite{BGP}, as follows
\begin{equation}\label{equazione_varn} \left\{
\begin{array}{ll}
dY^{n,\e}(t) = AY^{n,\e}(t)\,dt  +G\hat{B}^n(t) Y^{n,\e}(t)\,dt  + \sum_{i=1}^m G \hat{\Sigma}^{i,n}(t) Y^{n,\e }(t) d\, W ^i(t)  + \sum_{i=1}^m G \delta \Sigma ^{i}(t)d\, W ^i(t),\\ 
Y^{n,\e}(0)=0\in H,
\end{array}
\right.
\end{equation}
and \vspace{-0.2truecm}
\begin{equation} \label{equazione_varII_n}\left\{
\begin{array}{ll}
dZ^{n,\e}(t)  = [AZ^{n,\e}(t) +G \hat{ B}^n(t) Z^{n,\e}(t)  +    G\delta B (t)  + \frac{1}{2} G \tilde{B}^n(t) Y^{n,\e}(t)^2]\, dt  \\  + 
 G\sum_{i=1}^n [ \frac{1}{2} \tilde{\Sigma}^{i,n}(t)  Y^{n,\e}(t)^2 + \hat{\Sigma}^{i,n}(t) Z^{n,\e} (t)  + 
\delta^{i,n}\hat{\Sigma}^{i,n}(t)  Y^{n,\e}(t)]\, d W^i(t),  \quad t \in [0,T], \vspace{0.3cm} \\ 
 Z^{n,\e}_0 \,=\,0 \in H.
\end{array} \right.
\end{equation}
We underline again that these equations are not equations for the first variations themselves, but only approximations of the first variation equations of the original state equation.
Such equations take values in $H$ and have a unique solution, see theorem \ref{teo:var-stato}, that is given by \vspace{-0.2truecm}
\begin{equation*}
    Y^{n,\e}(\cdot)=\left( \begin{array}{l}
     y^{n,\e}(\cdot)\\
         y^{n,\e}(\cdot+\theta)_{\theta\in[-\mathtt{d},0]}
    \end{array}\right) , \,
    Z^{n,\e}(\cdot)=\left( \begin{array}{l}
     z^{n,\e}(\cdot)\\
       z^{n,\e}(\cdot+\theta)_{\theta\in[-\mathtt{d},0]}
    \end{array}\right)
\end{equation*}
where  $y^{n,\e}, z^{n,\e}$ are solutions of \eqref{yn}, \eqref{zn}.
Thus  \eqref{stima sup yn} and \eqref{conv y} can be rewritten in terms $Y^{n,\e}$, $Z^{n,\e}$, $Y^{\e}$ and $Z^{\e}$ as follows \vspace{-0.2truecm}
 \begin{equation*}
        \sup_{n\geq 1}  \E \sup_{0\leq t \leq T} |Y^{n,\e}(t)|^p < +\infty, \qquad
\sup_{n\geq 1}  \E \sup_{0\leq t \leq T} |Z^{n,\e}(t)|^p < +\infty
    \end{equation*}  
    Moreover it holds:\vspace{-0.2truecm}
    \begin{equation*}
        \lim_{n \to \infty} \E \sup_{0 \leq t \leq T} |Y^{n,\e}(t)- Y^{\e}(t)|^2=0, 
   \qquad
        \lim_{n \to \infty} \E \sup_{-\mathtt{d} \leq t \leq T} |Z^{n,\e}(t)- Z^{\e}(t)|^2=0, 
    \end{equation*}
    
\subsection{Duality relation II: second order adjoint equation}

In this Section we introduce the second order adjoint equation  of \eqref{equazione_varn}, in order to do this we consider the hamiltonian, see \eqref{hamiltoniana ridotta}, 
\begin{align}\label{hamiltoniana ridotta bis}\calh^a(t, p^a(t),q^a(t))=\sum_{j=1}^d b^j(t) p^{a}_j(t)
 + \sum_{j=1}^d \sum_{i=1}^m\sigma^ {j,i}(t) q^a_{i,j}(t) 
-\ell(t):=B^*(t) {p}^a(t) +\Sigma^*(t)q^a(t)-L(t)
\end{align}
Notice that $\mathcal H^a$ is evaluated at $(p^a,q^a)$, which is the solution of the ABSDE \eqref{anticipante-sec}. 
   We recall that  $\mathcal H:[0,T]\times H\times L (\R^m;H)\rightarrow \R$ while $\mathcal H^a:[0,T]\times \R^d\times L (\R^m;\R^d)\rightarrow \R$, so here we define $p(t)=\left(\begin{array}{c}
         p^a(t) \\0
        \end{array}\right), \,q(t)=\left(\begin{array}{c}
         q^a(t) \\0
        \end{array}\right)$ and $\calh(t, G^*p(t),G^*q(t))=\calh^a(t, p^a(t),q^a(t))$. Notice that the second order derivative of $\calh$ can be written also in terms of $( p^a,q^a)$:
   \begin{align} \label{hamiltoniana ridotta de}
 \mathcal H_{XX}(t)&= \mathcal{H}_{XX}(t, p(t), q(t))  =  
 B_{XX}^*(t) G^*p(t)  +  \sum_{i=1}^m \Sigma_{XX}^{i}(t)^* G^*q^{i}(t)  - L_{XX} (t) \\ \nonumber
&=B_{XX}^*(t) p^a(t)  +  \sum_{i=1}^m \Sigma_{XX}^{i}(t)^* q^{a}(t)  - L_{XX} (t).
\end{align}
    We also introduce a suitable, {\it symmetric} approximation for the second order derivative
 \begin{align}\label{Hn}\calh^{n}_{XX}(t)=(\tilde{B}^n(t))^* {p}^a(t) +(\tilde {\Sigma}^n(t))^*q^a(t)-\tilde{L}^n(t).
\end{align}  

\begin{proposition} \label{approxII}
For every $n \in \mathbb{N}$ we consider the BSDE
\begin{equation}\label{second_adjoint_n} \left\{
\begin{array}{ll} 
-dP^n (t)= \left( P^n(t)A+ A^*P^n(t) + P^n(t) G\hat{B}^n (t)+ (\hat{B}^n (t))^* G^*P^n(t)\right) \, dt +  \calh^n_{XX}(t) \, dt +\\ \\
\qquad \quad \ \  +\sum_{i
=1}^m [ 
 \hat\Sigma^{i,n}(t)^*  G^* P^n (t) G  \hat\Sigma^{i,n}(t)]\, dt \\  \\
+ \sum_{i=1}^m \left(  \hat\Sigma^{i,n}(t) ^*  G^* Q^{i,n} (t)  +  
 Q^{i,n}(t)G  \hat\Sigma^{i,n}(t)\right)\,d t  
+ \sum_{i=1} ^m Q^{i,n}(t)dW^i(t)
  \vspace{1mm}\\  
   P^n(T)=- H_{XX} (T).
\end{array}\right.
\end{equation}
  Then under Hypotheses \ref{ipotesi} and \ref{ipotesi-costo}
  there exists a unique mild solution $(P^n,Q^n)$ to equation \eqref{second_adjoint} such that 
$ P^n\in L^2_{\mathcal{F}}(\Omega; C([0,T]; S^2(H))) $  and $
Q^n \in L^2_{\mathcal{F}}(\Omega\times [0,T]; L(\R^m;S^2(H)))$.
\end{proposition}
\dim
The proof follows the lines of the proof of Theorem \ref{teo-exist-IIorder}.
\qed

\smallskip

Now we write $P^n(t)$ analogously to what done in \eqref{P-forma-matriociale}
\begin{equation*}
P^n(t)=\left(\begin{array}{cc}
     P^n_{0,0}(t)&P^n_{0,1} (t) \\
 P^n_{1,0}(t)&P^n_{1,1}(t)
\end{array}
\right) 
\end{equation*}
and we will prove the convergence of $P_{0,0}^n$ to a process $ P_{0,0}$. Although we are not able to characterize the limit of equation \eqref{second_adjoint_n}, this suffices to write down the stochastic maximum principle in its necessary form.

For this purpose we consider, for every $ h \in \R^d$ and $ s \in [0,T]$, the following equation
\begin{equation}\label{Yeq1}\left\{
\begin{array}{ll}
dY(t) = AY(t)\,dt  +G\hat{B}(t) Y(t)\,dt  + \sum_{i=1}^m G \hat{\Sigma}^{i}(t) Y(t) d\, W ^i(t),  \  t \in [s,T],\\ 
Y(s)=\left(\begin{array}{l}h\\0\end{array}\right)\in \R^d\times D([-\mathtt{d},0],\R^d),
\end{array}
\right.
\end{equation}
 Here $D([-\mathtt{d},0],\R^d)$ denotes the space of cadlag functions, indeed up to time $t=\mathtt{d} +s$ the past trajectory can have one discontinuity of cadlag type, see also \cite[Section 2]{FlZa}, for a discussion on the choice of the spaces. We refer again to \cite[Section 2]{FlZa} where it is discussed the fact that $A$ generates a semigroup of operators in $\R^d\times D([-\mathtt{d},0],\R^d)$. We don't really need to work with equation \eqref{Yeq1} reformulated in this space, we will use only that equation \eqref{Yeq1} has a solution $ Y^{s,h}$ that corresponds to
\begin{equation*}
    Y^{s,h}(\cdot)=\left( \begin{array}{l}
     y^{s,h}(\cdot)\\
         y^{s,h}(\cdot+\theta)_{\theta\in[-\mathtt{d},0]}
    \end{array}\right)\in\R^d\times D([-\mathtt{d},0],\R^d)
\end{equation*}
and $y^{s,h}$ solves 
\begin{equation}\label{y.dato}
\left\{ \begin{array}{ll}
d{y}(t) = \Big[\displaystyle b_x(t) {y}(t) + b_y(t) \int_{-\mathtt{d}}^ 0{y} (t+ \theta) \mu_b(d\theta) \Big] \, dt   \\ \\  \qquad\quad + \Big[\displaystyle\sigma_x(t) {y}(t) + \sigma_y(t) \int_{-\mathtt{d}}^0 {y} (t+ \theta) \mu_\sigma(d\theta) ] d W(t), \quad    t \in [s,T], \\ \\
  {y}(s)=h, \quad y(\theta)=0, \quad -d \leq \theta <s.
\end{array} \right.
\end{equation}

We need to introduce also its approximation $Y^{n,s,h}$, that  solves, for each $h \in \R^d$ and $n \geq 1$:
\begin{equation}\label{yn-bis-abstr}\left\{
\begin{array}{ll}
dY(t) = AY(t)\,dt  +G\hat{B}^n(t) Y(t)\,dt  + \sum_{i=1}^m G \hat{\Sigma}^{i,n}(t) Y(t) d\, W ^i(t),  \quad  t \in [s,T],\\ 
Y(s)=\left(\begin{array}{l}
h\\ 0
\end{array}\right)\in H,\end{array}
\right.
\end{equation}
and we consider the solutions of equations \eqref{Yeq1} and \eqref{yn-bis-abstr} defined on the whole $[0, T]$ by extending them to $[0,s)$ setting $Y(r)=\left(\begin{array}{l}
0\\ 0\end{array}\right),\, r\in [0,s)$.
\newline Notice that since we consider the coefficients $\hat B^n,\, \hat \Sigma^{i,n}$ with regularized measures, we have that equation \eqref{yn-bis-abstr} evolves in $H$, and also in $\R^d\times D([-\mathtt{d},0],\R^d)$ as discussed for equation \eqref{Yeq1}.
Moreover as usual we have
\begin{equation*}
    Y^{n,s,h}(\cdot)=\left( \begin{array}{l}
     y^{n,s,h}(\cdot)\\
         y^{n,s,h}(\cdot+\theta)_{\theta\in[-\mathtt{d},0]}
    \end{array}\right)\in\R^d\times D([-\mathtt{d},0],\R^d),
\end{equation*}
where   $y^{n,s,h}$ solves 
\begin{equation}\label{yn.dato}
\left\{ \begin{array}{ll}
d{y}(t) = \Big[\displaystyle b_x(t) {y}(t) + b_y(t) \int_{-\mathtt{d}}^ 0{y} (t+ \theta) f^n_b(\theta)d\theta \Big] \, dt   \\ \\  \qquad\quad + \Big[\displaystyle\sigma_x(t) {y}(t) + \sigma_y(t) \int_{-\mathtt{d}}^0 {y} (t+ \theta) f^n_b(\theta)d\theta  ] d W(t), \quad    t \in [s,T], \\ \\
  {y}(s)=h, \quad y(\theta)=0, \quad -d \leq \theta <s.
\end{array} \right.
\end{equation}
Notice that equations \eqref{y.dato} and \eqref{yn.dato} has an initial datum which is not continuous at $0$, so we cannot apply directly Theorem \ref{teo-mohammed}, but arguing as in the book \cite{Moha}, see also \cite{FlZa}, we can easily deduce that equations \eqref{y.dato} and \eqref{yn.dato} admit a unique cadlag solution, with the past trajectory that has at most one discontinuity point, and  the following estimate holds true for $p\geq 1$ where $x(\cdot)$ in turn denotes the solution to equation \eqref{y.dato} or \eqref{yn.dato}.
\begin{equation} \label{stima_dati_stato_cadlag}
\E \sup_{t \in [-\mathtt{d},T]} |x(t)|^p \leq c (1+ |h|^p)
\end{equation}

 The analogous results as in Proposition \ref{prop.conv.yn-zn}, follows: 
 
 \begin{proposition} \label{prop.conv yhn}
     Under Hypotheses \ref{ipotesi} and \ref{ipotesi-costo}, for every  $ h \in \R$ and $ s \in [0,T]$ and $ p \geq 1$,
      \begin{equation}\label{stima sup Ynh}
        \sup_{n\geq 1}\,  \E \sup_{0\leq t \leq T} |Y^{n,s,h}(t)|^p=  \sup_{n\geq 1} \, \E \sup_{-\mathtt{d}\leq \rho \leq T} |y^{n,s,h}(\rho)|^p < +\infty, 
    \end{equation}

    \begin{equation}\label{stima sup Yh}
        \sup_{n\geq 1}  \E \sup_{0\leq t \leq T} |Y^{s,h}(t)|^p =\sup_{n\geq 1}  \E \sup_{-\mathtt{d}\leq \rho \leq T} |y^{s,h}(\rho)|^p  < +\infty, 
    \end{equation}
     \begin{equation}\label{lim Y h} \lim_{n \to +\infty}  \E\sup_{t \in [s,T]} | Y^{n,s,h}(t) -  Y^{s,h}(t)|^ 2 = \lim_{n \to +\infty}  \E\sup_{-\mathtt{d}\leq \rho \leq T } | y^{n,s,h}(\rho) -  y^{s,h}(\rho)|^2=0. 
     \end{equation}
And:
\begin{equation} \label{conv debole b}
 \lim_{n \to \infty} \E \int_0^T \Big| \int_{-\mathtt{d}}^0
[y^{h,s} (t+\theta)f^n_b (\theta) \, d \theta -  y^{h,s} (t+\theta)\mu_b(d\theta)  \Big|^2 \, dt =0. 
\end{equation}
\begin{equation} \label{conv debole sigma}
 \lim_{n \to \infty} \E \int_0^T \Big| \int_{-\mathtt{d}}^0
[y^{h,s} (t+\theta)f^n_\sigma (\theta) \, d \theta -  y^{h,s} (t+\theta)\mu_\sigma(d\theta)  \Big| ^2 \, dt =0.
\end{equation}
 \end{proposition}
 \dim 
See the Appendix \ref{appendix}.
\qed
\begin{proposition}\label{Conv P00}
Under  hypotheses \ref{ipotesi} and \ref{ipotesi-costo}, for every $ h \in \R$ and for a.a.  $s\in[0,T]$.
\begin{align*}
 \lim_{n\to \infty}   \<P^n_{0,0}(s)h,h\> = \E^{\F_s}\<H_{XX}(T)Y^{s,h}(T), Y^{s,h}(T)\>-\E^{\F_s}\int_s^T \<\calh_{XX}(t)Y^{s,h}(t), Y^{s,h}(t)\>\,dt,
    \end{align*}
\end{proposition}
\dim Remember that, in this Section, we are considering a final cost $h$ depending only on $x(T)$, while we can still consider a current cost depending on the past trajectories. Since the dependence on the past trajectory on $\ell$ gives raise to terms that can be treated as the ones originating from the dependence on the past trajectory in $b$ and $\sigma$, for the sake of clearness, without affecting the generality of the result, we will suppose that $\ell$ doesn't depend of the past of the trajectory. Therefore $ \calh^{n}_{XX}$ reads as follows:
\begin{align}\label{Hnshort}\calh^{n}_{XX}(t)=(\tilde{B}^n(t))^* {p}^a(t) +(\tilde {\Sigma}^n(t))^*q^a(t)-L_{XX}(t).
\end{align}  
Passing through the Yosida approximants and
by applying the It\^o formula we have 
\begin{align*}
     d \<P^n(t)Y^{n,s,h}(t), Y^{n,s,h}(t)\>=\-\<\calh^n_{XX}Y^{n,s,h}(t), Y^{n,s,h}(t)\>, 
\end{align*}
where  $Y^{n,s,h}$ solves \eqref{yn-bis-abstr}. 
So integrating between $s$ and $T$ and taking the conditional expectation we get 
\begin{align} \label{rappresentazione n}
    \<P^n_{0,0}(s)h,h\>= \E^{\F_s}\<H_{XX}(T)Y^{n,s,h}(T), Y^{n,s,h}(T)\>-\E^{\F_s}\int_s^T \<\calh^n_{XX}(t)Y^{n,s,h}(t), Y^{n,s,h}(t)\>\,dt.
\end{align}
In order to let $n \to +\infty $ in \eqref{rappresentazione n} we notice that:
\begin{align*}
& \<\calh^n_{XX}(t)Y^{n,s,h}(t), Y^{n,s,h}(t)\> - \<\calh_{XX}(t)Y^{s,h}(t), Y^{s,h}(t)\>=  \<[\calh^n_{XX}(t) - \calh_{XX}(t)]Y^{s,h}(t), Y^{n,s,h}(t)\>  \\&  +\<\calh^n_{XX}(t)[Y^{n,s,h}(t) -Y^{s,h}(t)], Y^{n,s,h}(t) \>   +   \<\calh_{XX}(t)Y^{s,h}(t), Y^{n,s,h}(t) -Y^{s,h}(t) \>\   \\&  = I_1+I_2+I_3
\end{align*}
Thanks to our hypotheses, in particular \eqref{stima sup yn} and  the fact that $|\mu_n| \leq K$, we have that:
\begin{align*}
    I_1 &\leq |[\calh^n_{XX}(t) - \calh_{XX}(t)]Y^{s,h}(t)| | Y^{n,s,h}(t)| \\ &\leq c \sup_{\rho \in [-\mathtt{d},T]}  |  y^{n,s,h}(\rho)||p^a(t)|  \Big|\int_{-\mathtt{d}}^0
[y^{h,s} (t+\theta)f^n_b (\theta) \, d \theta -  y^{h,s} (t+\theta)\mu_b(d\theta)   \Big| \\& +  c \sup_{\rho \in [-\mathtt{d},T]}  |  y^{n,s,h}(\rho)| |q^a(t)|  \Big|\int_{-\mathtt{d}}^0
[y^{h,s} (t+\theta)f^n_\sigma (\theta) \, d \theta -  y^{h,s} (t+\theta)\mu_\sigma(d\theta)   \Big|
\end{align*}
Thus we deduce that
\begin{align*}
 &\E\int_s^T| \<[\calh^n_{XX}(t) - \calh_{XX}(t)]Y^{s,h}(t), Y^{n,s,h}(t)\>|\,dt \nonumber\\ &\leq c \,  \Big \{ \E \Big[ \Big( \int_s^T   |p^a(t)|^2 \, dt \Big)^{1/2}   \sup_{\rho \in [-\mathtt{d},T]}  |  y^{n,s,h}(\rho)|\Big( \int_s^T \Big|\int_{-\mathtt{d}}^0
[y^{h,s} (t+\theta)f^n_b (\theta) \, d \theta -  y^{h,s} (t+\theta)\mu_b(d\theta)   \Big|^2 \, dt \Big)^ {1/2}\Big] \\& 
+   \E \Big[ \Big( \int_s^T  |q^a(t)|^2 \, dt \Big)^{1/2}  \sup_{\rho \in [-\mathtt{d},T]}  |  y^{n,s,h}(\rho)|\Big( \int_s^T \Big|\int_{-\mathtt{d}}^0
[y^{h,s} (t+\theta)f^n_{\sigma} (\theta) \, d \theta -  y^{h,s} (t+\theta)\mu_{\sigma}(d\theta)   \Big|^2 \, dt \Big)^ {1/2}\Big] \Big \} 
\\& \leq   c \Big( \E \int_0^T   |p^a(t)|^2  +  |q^a(t)|^2  \, dt \Big)^ {1/2}  \Big( \E \int_0^T  \sup_{\rho \in [-\mathtt{d},T]}  |  y^{n,s,h}(\rho)|^2  
 \times\\&  \Big[ \Big|\int_{-\mathtt{d}}^0
y^{h,s} (t+\theta)f^n_b (\theta) \, d \theta -  y^{h,s} (t+\theta)\mu_b(d\theta)   \Big|^2 +\Big|\int_{-\mathtt{d}}^0
[y^{h,s} (t+\theta)f^n_\sigma (\theta) \, d \theta -  y^{h,s} (t+\theta)\mu_\sigma(d\theta)   \Big|^2 \Big]\, dt \Big)^ {1/2}
\\ &  \leq  c \Big( \E \int_0^T   |p^a(t)|^2  +  |q^a(t)|^2  \, dt \Big)^ {1/2} 
\Big( \E  \sup_{\rho \in [-\mathtt{d},T]}  |  y^{n,s,h}(\rho)|^8  \Big) ^{1/8} \Big( \E \sup_{\rho \in [-\mathtt{d},T]}  |  y^{s,h}(\rho)|^4\Big) ^{1/8} \times \\&
\Big( \E   \int_0^T \Big|\int_{-\mathtt{d}}^0
[y^{h,s} (t+\theta)f^n_b (\theta) \, d \theta -  y^{h,s} (t+\theta)\mu_b(d\theta)   \Big| ^2 +\Big|\int_{-\mathtt{d}}^0
[y^{h,s} (t+\theta)f^n_\sigma (\theta) \, d \theta -  y^{h,s} (t+\theta)\mu_\sigma(d\theta)   \Big| ^2\, dt \Big)^ {1/4}
\end{align*}
where the constant $c$ that may change from line to line but remains always independent of $n$, hence by  Proposition \ref{prop.conv yhn}, we  deduce that \vspace{-0.2truecm}
\begin{equation}
 \lim_{n \to +\infty}\E\int_s^T \<[\calh^n_{XX}(t) - \calh_{XX}(t)]Y^{s,h}(t), Y^{n,s,h}(t)|\>\,dt=0
\end{equation}
Thus, again by Proposition \ref{prop.conv yhn} and Dominated Convergence, we deduce that:\vspace{-0.2truecm}
\begin{equation}
 \lim_{n \to +\infty}\E^{\F_s}\int_s^T \<\calh^n_{XX}(t)Y^{n,s,h}(t), Y^{n,s,h}(t)\>\,dt=
 \E^{\F_s}\int_s^T \<\calh_{XX}(t)Y^{s,h}(t), Y^{s,h}(t)\>\,dt
\end{equation}
While for $I_2$, and similarily for $I_3$, we notice that\vspace{-0.2truecm} \begin{equation}
    I_2 \leq c  [|p^a(t)|  +  |q^a(t)| +1]   \sup_{\rho \in [-\mathtt{d},T]} | y^{n,h,s} (\rho) - y^{h,s} (\rho) |  \sup_{\rho \in [-\mathtt{d},T]} | y^{n,h,s} (\rho)|
\end{equation}
Therefore \vspace{-0.6truecm}
\begin{align*}
 &\E\int_s^T|  \<\calh^n_{XX}(t)[Y^{n,s,h}(t) -Y^{s,h}(t)], Y^{n,s,h}(t) \> | \, dt  \\\leq c  & 
\Big( \E\int_s^T  [|p^a(t)| ^2 +  |q^a(t)|^2  +1] \,dt \Big)^{1/2}
 \Big( \E \sup_{\rho \in [-\mathtt{d},T]}  |  y^{s,h}(\rho)|^2\Big) ^{1/2} \times \\&
  \Big(\E \sup_{\rho \in [-\mathtt{d},T]} | y^{n,h,s} (\rho) -y^{h,s} (\rho) |^2\Big)^ {1/2}
 \end{align*}
Thus by  Proposition \ref{prop.conv yhn},  in view of \eqref{lim Y h}, we deduce that:\vspace{-0.2truecm}
\begin{equation}
 \lim_{n \to +\infty}\int_s^T  \<\calh^n_{XX}(t)[Y^{n,s,h}(t) -Y^{s,h}(t)], Y^{n,s,h}(t) \> \,dt=
 0.
\end{equation}
We eventually let $n$ tend to $\infty$ in \eqref{rappresentazione n} and get for every $s \in [0,T]$, $\forall h \in H$:\vspace{-0.2truecm}
 \begin{align}\label{eq:P00-1}
    \lim_{n \to +\infty}\< P^n_{0,0}(s)h,h\>= \E^{\F_s}\<H_{XX}(T)Y^{s,h}(T), Y^{s,h}(T)\>-\E^{\F_s}\int_s^T \<\calh_{XX}(t)Y^{s,h}(t), Y^{s,h}(t)\>\,dt, \ \P-a.s..
\end{align}
\qed

We then deduce the characterization  for the limit.
\begin{theorem}\label{teoP00}
    Under hypotheses   \ref{ipotesi} and \ref{ipotesi-costo}, we have that there exists a unique couple $P_{0,0} \in L^2_\F (\Omega \times [0,T]; \R^{d\times d})$   and $Q_{0,0}(\cdot, \cdot) \in L^2([0,T];  L^2_\F (\Omega \times [0,T]; \R^{d\times d}))$ such that, for every $s \in [0,T]$:
    \begin{equation}
     \< P_{0,0}(s)h,h\> =   \<H_{XX}(T)Y^{s,h}(T), Y^{s,h}(T)\>-\int_s^T \<\calh_{XX}(t)Y^{s,h}(t), Y^{s,h}(t)\>\,dt - \int_s^T Q_{0,0}(s,t) \, dW(t)\vspace{-0.2truecm}
    \end{equation}
    and, for every $ s \in [0,T]$\vspace{-0.2truecm}
    \begin{equation}
        P_{0,0} (s) = \lim_{n\to \infty} P^n_{0,0}(s), \qquad \qquad \P -a.s.
    \end{equation}
\end{theorem}
\dim
 Let us define a matrix valued process as follows:\vspace{-0.2truecm}
$$\<\phi(s)h,k\> = : \<H_{XX}(T)Y^{s,h}(T), Y^{s,k}(T)\>-\int_s^T \<\calh_{XX}(t)Y^{s,h}(t), Y^{s,k}(t)\>\,dt, \qquad h,k \in \R^d\vspace{-0.2truecm}$$ 
 We  easily  deduce,  thanks to hypotheses \ref{ipotesi} and \ref{ipotesi-costo} and to the regularity of $Y^{s,h}$, that  $ \phi \in L^2([0,T]\times \Omega; \R^{d \times d}).$  Thus, by the Martingale Representation Theorem, see also \cite [Theorem 2.2]{Yong 2006},  there exists  $(P_{0,0} (\cdot), Q_{0,0}(\cdot,\cdot)) \in   L^2_\F (\Omega \times [0,T]; \R^{d\times d})\times  L^2([0,T];  L^2_\F (\Omega \times [0,T]; L(\R^m; \R^{d\times d}))) $, such that, for all $s \in [0,T]$:\vspace{-0.4truecm}
\begin{equation}
      P_{0,0}(s) =   \phi(s) - \int_s^T Q_{0,0}(s,t) \, dW(t)
    \end{equation}
Moreover, 
the following holds, for every $h \in \R^d$:\vspace{-0.2truecm}
\begin{equation}\label{teoP00formula}
      \<P_{0,0}(s)h,h\> =  \E ^{\F_s} \<H_{XX}(T)Y^{s,h}(T), Y^{s,h}(T)\>- \E ^{\F_s} \int_s^T \<\calh_{XX}(t)Y^{s,h}(t), Y^{s,h}(t)\>\,dt 
    \end{equation}
Thus, by \eqref{eq:P00-1}, we get that $P_{0,0}(s) =  \lim_{n\to \infty} P^n_{0,0}(s)$.
\qed\vspace{-0.2truecm}
\subsection{Maximum principle}\label{max-princ-gen}


Now we are able to prove a version of the Stochastic Maximum Principle, in its necessary form, for the control problem with state equation and cost functional given by \eqref{equazione_stato} and \eqref{costo}, respectively and with measures $\mu_b$ and $\mu_\sigma$ that
are general finite regular measures on $[-\mathtt{d},0]$.  
\noindent We consider $y^\e,\,z^\e$ which are respectively the first and second order variation of the stochastic delay equation \eqref{equazione_stato} and we set\vspace{-0.2truecm}
\begin{equation*}
    Y^\e(\cdot)=\left( \begin{array}{l}
     \quad y^\e(\cdot)\\
         y^\e(\cdot+\theta)_{\theta\in[-\mathtt{d},0]}
    \end{array}\right),\,
    Z^\e(\cdot)=\left( \begin{array}{l}
     \quad z^\e(\cdot)\\
       z^\e(\cdot+\theta)_{\theta\in[-\mathtt{d},0]} 
    \end{array}\right) \in\R^d\times C([-\mathtt{d},0],\R^d).\quad 
\end{equation*}
We will consider the expansion of the cost  (\ref{costo-espansione-I-gen}) given in Corollary \ref{espansioneI-gen}.
Taking into account that\vspace{-0.2truecm}
\begin{equation*}
    Y^{n,\e}(\cdot)=\left( \begin{array}{l}
     y^{n,\e}(\cdot)\\
         y^{n,\e}(\cdot+\theta)_{\theta\in[-\mathtt{d},0]}
    \end{array}\right) , \,
    Z^{n,\e}(\cdot)=\left( \begin{array}{l}
     z^{n,\e}(\cdot)\\
       z^{n,\e}(\cdot+\theta)_{\theta\in[-\mathtt{d},0]}
    \end{array}\right)\in\R^d\times C([-\mathtt{d},0],\R^d),\quad 
\end{equation*}
  we have 
  
 \begin{align}\label{costo-espansione-gen-n}
J(u)&-J(u^\e)=\E\int_0^T \left(-
\delta \ell (t)+ \<p^a(t),\delta b(t) +\<q^a(t),\delta \sigma(t) \>\>\right) I_{E_\e}(t)\,dt \\ 
 &\nonumber-\frac{1}{2} \E \<h_{xx}(X(T))y^\e(t),y^\e(t)\> 
-\frac{1}{2}\E \int_0^T \<\ell_{xx}(t)y^\e(t),y^\e(t)\>\, dt \\ \nonumber& 
 -\frac{1}{2}\E \int_0^T \int_{\mathtt{d}}^ 0\int_{\mathtt{d}}^ 0\ell_{yy}(t) y^\e(t+\theta_1)y^\e(t+\theta_2){\mu} _\ell(d\theta_1) {\mu} _\ell(d\theta_2)\,dt
-\E\dis\int_0^T\int_{-\mathtt{d}}^0 \ell_{xy}(t)y^\e(t)  y^{\e}(t+\theta)\mu_\ell(d\theta)\,dt\\&
\nonumber
+\E\dis\int_0^T\frac{1}{2} \<p^a(t), b_{xx}(t)y^\e(t)y^\e(t)+\int_{-\mathtt{d}}^ 0\int_{-\mathtt{d}}^ 0b_{yy}(t) y^\e(t+\theta_1)y^\e(t+\theta_2){\mu} _b(d\theta_1) {\mu} _b(d\theta_2)\>\,ds\\  \nonumber
&\quad +\E\dis\int_0^T\<p^a(t),\int_{-\mathtt{d}}^0 b_{xy}(t)y^\e(t)  y^{\e}(t+\theta)\mu_b(d\theta)\>\,dt\\  \nonumber
 &\quad+ \E\dis\int_0^T\frac{1}{2}\<q^a(t), \sigma_{xx}(t)y^\e(t)y^\e(t)+\int_{-\mathtt{d}}^ 0\int_{-\mathtt{d}}^ 0\sigma_{yy}(t) y^\e(t+\theta_1)y^\e(t+\theta_2){\mu} _\sigma(d\theta_1) {\mu} _\sigma(d\theta_2)\>\,dt\\  \nonumber
&\quad +\E\dis\int_0^T\<q^a(t),\int_{-\mathtt{d}}^0 \sigma_{xy}(t)y^\e(t)  y^{\e}(t+\theta)\mu_\sigma(d\theta)\>\,dt
 +o(\e)\\ \nonumber
 &=\E\int_0^T \left(-
\delta \ell (t)+ \<p^a(t),\delta b(t) +\<q^a(t),\delta \sigma(t) \>\>\right) I_{E_\e}(t)\,dt \\ \nonumber
& -\frac{1}{2}\E 
 H_{XX} (T)(Y^{\e}(T))^2+
\frac{1}{2}\E \int_0^T \calh_{XX}(t)(Y^{\e}(t))^2\, dt+o(\e).
\end{align}
Now, by the convergences stated in Proposition \ref{prop.conv.yn-zn} and arguing exactly as in  Prop. \ref{Conv P00}, we get\vspace{-0.2truecm}
\begin{align*}
    &\E  \,
 H_{XX} (T)(Y^{\e}(T))^2 - \E \int_0^T \calh_{XX}(t)(Y^{\e}(t))^2\,dt = \\
& \nonumber \lim_{n\to+\infty}\left[\E  \, 
 H_{XX} (T)(Y^{n,\e}(T))^2-
\E \int_0^T \calh_{XX} ^n(t)(Y^{n,\e}(t))^2\, dt \right].
\end{align*}
In the following Theorem the stochastic maximum principle is written in terms of the Hamiltonian function $\calh$ defined in \eqref{hamiltoniana ridotta}, see also \eqref{hamiltoniana ridotta bis}.
\begin{theorem}\label{teoMaxPrinc-gen}
Let Assumptions (A.1) and (A.2) be satisfied and suppose that $(X,u)$ is an optimal pair for the control problem, 
and let us consider $u^\e$ defined in (\ref{spike_def}).
There exist two pairs of processes $ (p^a,q^a) \in L^2_{\mathcal{F}}(\Omega; C([0,T];\R^d)) \times
 L^2_{\mathcal{F}}(\Omega\times [0,T];L(\R^m,\R^d))$, solution to the first order adjoint equation \eqref{eq:adjoint-I}, 
such that the following variational inequality holds
 \begin{align*}
  &\calh(t, X(t),  u(t),p^a(t),q^a(t))-   \calh(t, X(t), v,p^a(t),q^a(t)) \\ &+\frac{1}{2} \sum_{i=1}^m
  [(\delta \Sigma^i(t))^* P_{0,0}(t) \delta \Sigma^i(t)]\leq 0
  \quad \forall\,v\in U,\quad \text{a.e. } t\in [0,T],\, \P-\text{a.s.}.\vspace{-0.2truecm}
 \end{align*}
 Equivalently, for $ \text{ a.e. } t\in [0,T],\, \P-\text{a.s.}$,\vspace{-0.2truecm}
 \begin{align*}
  &\calh(t, X(t),  u(t),p^a(t),q^a(t)) +\frac{1}{2} \sum_{i=1}^m [(\delta \Sigma^i(t, X(t),u(t)))^* P_{0,0}(t) \delta \Sigma^i(t, X(t),u(t))]\\ 
  & \;=\inf_{v\in U}  \calh(t, X(t),v,p^a(t),q^a(t))
  +\frac{1}{2}  \sum_{i=1}^m[(\Sigma^i(t, X(t),v)^*P_{0,0}(t)\Sigma^i(t, X(t),v)] .\end{align*}
\end{theorem}
\dim Starting from \eqref{costo-espansione-gen-n} we can write, taking into account \eqref{relaggiunta II}, and remark \ref{rm:prima-componente}

\begin{align}\label{costo-espansione-II-n}
 J(u)-J(u^\e) 
 &=\E\int_0^T \left(-
\delta \ell (t)+ \<p^a(t),\delta b(t) \>  +\<q^a(t),\delta \sigma(t) \>\right) I_{E_\e}(t)\,dt\\ \nonumber&+ \frac{1}{2}\lim_{n\to+\infty}  \E \int_0^T\left[
  \sum_{i=1}^m  \delta \Sigma ^{i} (t) ^*P^n_{0,0}(t)\Sigma ^{i} (t)) I_{E_\e}(t)\, dt+o(\e)\right]\\ \nonumber&=
 \E\int_0^T \big[\calh(t, X(t),  u(t),p^a(t),q^a(t))-   \calh(t, X(t), v,p^a(t),q^a(t)) \\ \nonumber &+\frac{1}{2} \E \int_0^T\sum_{j=1}^m [(\delta \Sigma^j(t))^* P_{0,0}(t) \delta \Sigma^j(t)]\big] I_{E_\e(t)} \, dt +o(\e)
\end{align}
where in the last passage $P_{0,0}$ is completely identified thanks to theorem \ref{teoP00}.
Notice that the right hand side is not equal to $0$ only for $t\in E_\e$. So 
if we divide both sides by $\e$ and we let $\e\rightarrow 0$ we get \vspace{-0.2truecm}
\[
 0\geq \calh(t, X(t),  u(t),p^a(t),q^a(t))-   \calh(t, X(t), v,p^a(t),q^a(t)) +\frac{1}{2} \sum_{j=1}^m [(\delta \Sigma^j(t))^* P_{0,0}(t)\delta \Sigma^j(t)]\vspace{-0.2truecm}
\]
$\forall\,v\in U$, a.e. $t\in [0,T],\, \P-$ a.s., and this concludes the proof. 
\qed




   \appendix

   \section{An optimal portfolio problem}
\label{sec:opt-port}
We consider a generalized Black and Scholes market with one risky asset, whose price at time $t$ is denoted by $S(t)$ and whose past trajectory from time $t-d$ up to time $t$ is denoted by $S_t$, and one non-risky asset, whose price at time $t$ is
denoted by $B(t)$. The result can be extended to the case of a Black and Scholes market with $j$ risky assets, whose prices at time t are denoted by $S^i(t),i = 1, . . . , j$, and one non-risky asset: for the sake  of simplicity we limit here to the case of only one risky asset.
\newline  The evolution of the prices is given by the following stochastic delay differential equation in a complete probability space $(\Omega, \Fcal, \mP)$ : \vspace{-0.2truecm}
\begin{equation}\label{eq:price_evolution}
\left\lbrace
\begin{array}{l}
dS(t)=S(t )\left[ b(t, S(t),\displaystyle\int_{-\mathtt{d}}^0S(t+\theta )\mu_{ b}(d\theta)) ) dt+\sigma(t, S(t),\displaystyle\int_{-\mathtt{d}}^0S(t+\theta )\mu_{\sigma}(d\theta)) ) dW_t  \right],\vspace{0.1cm}\\\
S(0)=S_0,\,S(\theta)=\nu_0(\theta),\vspace{0.1cm}\\
dB(t)=r(t, S(t),\displaystyle\int_{-\mathtt{d}}^0S(t+\theta )\mu_{r}(d \theta))B(t) dt,\vspace{0.1cm}\\
B(0)=B_0
\end{array}
\right.
\end{equation}
where $W(t)$ is a standard Brownian motion in $\R$, $(\Fcal_t)_{t\geq 0}$ is the filtration generated by $W$ and augmented with null probability sets
The measures $\mu_{ b},\,\mu_{ \sigma},\,\mu_{ r}$ respectively in the drift $b$, in the diffusion $\sigma$ and in the rate $r$ are
given by
finite regular measures on $[-\mathtt{d},0]$. Notice that the presence of a regular measure allows to modulate the contribution of the past in the evolution of the price, for example giving more weight to the most recent values, see also \cite[Remark 2.3]{BGP} where   the advantage  of having general dependence on the past in the coefficients is discussed. 
\begin{Hypothesis}
\label{ip: b_r_sigma}
On $b, \,\sigma,\,r$ we make the following assumptions: 
\begin{itemize}
\item[i)]   $\mu_{ b}$ is a regular measure and $ b : [0,T]\times \R \times \R \rightarrow \R$  is measurable and 
for all $t\in[0,T]$ $ b(t,\cdot,\cdot)$ is twice differentiable with bounded derivatives;
\item[ii)] $\mu_{\sigma}$ is a regular measure and $\sigma : [0,T]\times \R\times \R\rightarrow \R$ is measurable 
and for all $t\in[0,T]$ $ \sigma(t,\cdot,\cdot)$ is twice differentiable with bounded derivatives;
\item[iii)] $
r: [0,T]\times \R\times \R\rightarrow \R$ is measurable
and for all $t\in[0,T]$ $ r(t,\cdot,\cdot)$ is twice differentiable with bounded derivatives.
\end{itemize}\vspace{-0.2truecm}
\end{Hypothesis}
Let us denote by $V(t)$ the value at time $t$ of the associated self-financing portfolio. We consider an optimal portfolio problem where we also allow consumption, and also the investors are allowed to take money from the portfolio $V$: in the model this is represented by a further control $c$, taking values in a non convex subset $\mathcal C \subset \R^+$, e.g. this is the case if $\mathcal C$ is discrete and its elements are multiples of a given quantity.

\smallskip

\noindent The state equation for the optimal portfolio is given by (see \cite{ElKa1997})
\begin{equation}\label{eq:portfolio}
\left\lbrace\begin{array}{ll}dV(t)= r (t,S(t), S_t))(V(t)- \pi(t)) dt-c(t)dt+\pi(t)\left[ b(t,S(t), S_t ) dt+\sigma(t,S(t), S_t ) dW_t )\right]
\\
V(0)=V_0,\,V(\theta)=\eta(\theta)\, \theta \in [-d,0),
\end{array}
\right.
\end{equation}
where for brevity we have written $f(t,S(t),S_t)=f(t,S(t),\int_{-\mathtt{d}}^0S(t+\theta )\mu_{f}(d\theta))$ for $f=b,\,\sigma,\,r$, and $S_t(\theta):=S(t+\theta),\, \theta \in [-\mathtt{d},0]$ denotes the past trajectory
 We will consider as admissible strategies  square-integrable, predictable investment $\pi\in {L}^2_\mathcal{F}(\Omega\times [0, T ] , \R)$.

The aim is to maximize the utility functional over the set of the admissible strategies; the utility is given by \vspace{-0.2truecm}
\begin{equation}\label{utility}
U(c)={\mathbb E} \int_0^T\left[U_1\left(t,V(t),c(t)\right) \right]\,dt +{\mathbb E} \left[U_2\left(V( T) \right) \right],
\end{equation}
where $U_1:[0,T]\times \R\rightarrow \R$ represents the utility from consumption and it is assumed to be Lipschitz continuous and twice differentiable in the second variable and  $U_2:\R\rightarrow \R$ represents the utility from the wealth at time $T$ and it is assumed to be Lipschitz continuous and twice differentiable. 

\noindent At any time $t\in [-\mathtt{d},T]$, the state $x(t)\in \R^2 $ is given by the pair
$
x(t)=\left(
\begin{array}{l}
S(t)\\
V(t)
\end{array}
\right).$
So the equation for $x$ is given by\vspace{-0.2truecm}
\begin{equation}\label{eq:pair_X}
\left\lbrace
\begin{array}{l}
d\left(
\begin{array}{l}
S(t)\\
V(t)
\end{array}
\right)=\left(\begin{array}{l}S(t ) b(t, S(t),S_t )\vspace{2mm}\\
 r (t,  S(t),S_t)(V(t)- \pi(t))-c(t)+\pi(t) b(t, S(t), S_t )\end{array}\right) dt\vspace{3mm}\\\qquad\,
 +\left(\begin{array}{l}S(t)\sigma(t,  S(t),S_t ) \vspace{2mm}\\
\pi(t) \sigma(t, S(t),  S_t )
\end{array}\right)dW_t \vspace{3mm}\\
\left(
\begin{array}{l}
S(0)\\
V(0)
\end{array}
\right)=\left(\begin{array}{l}
S_0\\
V_0
\end{array}
\right),\;
\left(
\begin{array}{l}
S(\theta)\\
V(\theta)
\end{array}
\right)=\left(\begin{array}{l}
\nu_0(\theta)\\
\eta(\theta)
\end{array}
\right)\, \theta \in [-d,0),
\end{array}
\right.
\end{equation}
and it turns out to be a controlled stochastic equation with delay in the state. Notice that we are considering the case of a general measure so we are in the framework of Section \ref{sec-gen-meas}, consequently as first adjoint equation we will consider an ABSDE, and we will characterize $P$ by means of Theorem \ref{teoP00}, namely by \eqref{teoP00formula}.
\newline Coming into the details, the first order adjoint processes are given by a pair of processes\vspace{-0.2truecm} $$(p^a,q^a)=\left(\left(\begin{array}{l}p^{a,1}\\p^{a,2}\end{array}\right),\left(\begin{array}{l}q^{a,1}\\q^{a,2}\end{array}\right)\right)\in L^2_{\mathcal{F}}(\Omega\times [0,T], \R^2)\times 
L^2_{\mathcal{F}}(\Omega\times [0,T], \R^2)$$ solution of the ABSDEs we are going to write, and that it turns out that the pair $(p^{a,1},q^{a,1})$ is identically $0$ since the cost functional doesn't depend on $x$, see also \cite[Section5]{GuaMas}.
The pair of processes $(p^{a,2},q^{a,2})\in \Lcal^2_{\mathcal{F}}(\Omega\times [0,T], \R)\times 
\Lcal^2_{\mathcal{F}}(\Omega\times [0,T], \R)$ satisfies the following equation :
 \begin{equation}\label{ABSDE-aggiuntaI_execution delay2} \left\lbrace\begin{array}{l}
p^{a,2}(t)=  \dis\int_t^Tp^{a,2}(s) \left(V(s)-\pi(s)\right)  r_x(s, S(s),S_{s})\mu_{r}(d\theta)\, ds\\ +\dis\int_t^T\E^{\Fcal_s}\dis\int_{-\mathtt{d}}^0 p^{a,2}(s-\theta) \left(V(s-\theta)-\pi(s-\theta)\right)  r_y(s-\theta, S(s-\theta),S_{s-\theta})\mu_{r}(d\theta)\, ds\\+\dis\int_t^T p^{a,2}(s)\pi(s) b_x(s,S(s), S_s)\, ds
+\dis\int_t^T\E^{\Fcal_s}\dis\int_{-\mathtt{d}}^0 p^{a,2}(s-\theta)\pi(s-\theta)  b_y(s-\theta, S(s-\theta),S_{s-\theta})\mu_{b}(d\theta)\, ds\\+\dis\int_t^T p^{a,2}(s)\pi(s) \sigma_x(s,S(s), S_s)\, ds
+\dis\int_t^T\E^{\Fcal_s}\dis\int_{-\mathtt{d}}^0 q^{a,2}(s-\theta) \pi(s-\theta)  \sigma_y(s-\theta,S(s-\theta), S_{s-\theta})\mu_{ \sigma}(d\theta)\, ds\\
+ \dis\int_t^T(U_{1})_v(s)ds+\dis\int_t^Tq^{a,2}(s)dW_s+ (U_2)_v\left(V( T)\right) \\p^{a,2}(T-\theta)=0, \; q^{a,2}(T-\theta)=0 \quad \forall \,\theta \in [-\mathtt{d},0).
 \end{array}\vspace{-0.2truecm}
 \right.
\end{equation}
In order to characterize $P_{0,0}$ we need to write down the equation for $y^{s,h}$, for $h=\left(
\begin{array}{l}
h^S\\
h^V
\end{array}
\right)\in \mathbb R^2$. 

\noindent We set $y^{s,h}(t)\equiv 0,$ tor $t<s$ while for $t\geq s$ $y^{s,h}(t):=\left(
\begin{array}{l}
y^{s,h,S}(t)\\
y^{s,h,V}(t)
\end{array}
\right)$ satisfies

\begin{equation}\label{eq:pair_ysh}
\left\lbrace
\begin{array}{l}
d\left(
\begin{array}{l}
y^{s,h,S}(t)\\
y^{s,h,V}(t)
\end{array}
\right)=\left(\begin{array}{l} S(t)\left[b_x(t, S(t),S_t )y^{s,h,S}(t)+b_y(t, S(t),S_t )\displaystyle\int_{-\mathtt d}^0y^{s,h,S}(t+\theta)\mu_b(d\theta)\right]\vspace{2mm}\\
 (V(t)- \pi(t))\left[r_x(t, S(t),S_t )y^{s,h,S}(t)+b_y(t, S(t),S_t )\displaystyle\int_{-\mathtt d}^0y^{s,h,S}(t+\theta)\mu_r(d\theta)\right]\end{array}\right) dt\vspace{3mm}\\\quad
 +\left(\begin{array}{l} b(t, S(t),S_t )y^{s,h,S}(t)\vspace{2mm}\\
 r (t,  S(t),S_t)y^{s,h,V}(t)+\pi(t)\left[b_x(t, S(t),S_t )y^{s,h,S}(t)+b_y(t, S(t),S_t )\displaystyle\int_{-\mathtt d}^0y^{s,h,S}(t+\theta)\mu_b(d\theta)\right]\end{array}\right) dt\vspace{3mm}\\\quad
 +\left(\begin{array}{l}\sigma(t, S(t),S_t )y^{s,h,S}(t)+S(t)\left[\sigma_x(t, S(t),S_t )y^{s,h,S}(t)+\sigma_y(t, S(t),S_t )\displaystyle\int_{-\mathtt d}^0y^{s,h,S}(t+\theta)\mu_\sigma(d\theta)\right]\vspace{2mm}\\
\pi(t) \left[\sigma_x(t, S(t),S_t )y^{s,h,S}(t)+\sigma_y(t, S(t),S_t )\displaystyle\int_{-\mathtt d}^0y^{s,h,S}(t+\theta)\mu_\sigma(d\theta)\right]
\end{array}\right)dW_t \vspace{3mm}\\
\left(
\begin{array}{l}
y^{s,h,S}(s)\\
y^{s,h,V}(s)
\end{array}
\right)=\left(\begin{array}{l}
h^S\\
h^V
\end{array}
\right),\;
\left(
\begin{array}{l}
y^{s,h,S}(s+\theta)\\
y^{s,h,V}(s+\theta)
\end{array}
\right)
=\left(\begin{array}{l}
0\\
0
\end{array}
\right),\, \theta \in [-d,0),
\end{array}
\right.
\end{equation}
Notice that here and in the following, for any function $f=f(t,S(t),S_t, V(t))$ ($f$ may also depends only on one, or in general not all, the arguments) we denote respectively by $f_x $ the derivative with respect to $S(t)$, by $f_y $ the derivative with respect to $S_t$, and by $f_v $ the derivative with respect to $V(t)$.

Next we also write down the Hamiltonian function:
\begin{align}\label{hamiltonian-example}
    \calh^a(t,x(t), \pi(t),c(t),p^a(t),q^a(t))&=\left[r (t,  S(t),S_t)(V(t)- \pi(t))-c(t)\right.\\ \nonumber
    &\left.+\pi(t) b(t, S(t), S_t )\right]p^{a,2}(t)+\pi(t)\sigma(t,  S(t),S_t )q^{a,2}(t)-U_1(t,V(t),c(t)).
\end{align}
Following Theorem \ref{teoP00}, namely \eqref{teoP00formula}, we are able to characterize $P_{0,0}$: for any $h=\left(
\begin{array}{l}
h^S\\
h^V
\end{array}
\right)\in \mathbb R^2$
\begin{align}\label{teoP00formula-example}
      \<P_{0,0}(s)h&,h\>
      =  \E ^{\F_s} \<(U_{2})_{vv}(T)y^{s,h,V}(T), y^{s,h,V}(T)\>- \E ^{\F_s} \int_s^T \<\pi(t) \sigma_{xx}(t, S(t), S_t )q^{a,2}(t)y^{s,h,S}(t), y^{s,h,S}(t)\>\,dt \\ \nonumber &- \E ^{\F_s} \int_s^T \<\left[r_{xx}(t,  S(t),S_t)(V(t)- \pi(t))+\pi(t) b_{xx}(t, S(t), S_t )\right]p^{a,2}(t)y^{s,h,S}(t), y^{s,h,S}(t)\>\,dt \\\nonumber &- \E ^{\F_s} \int_s^T\int_{-\mathtt d}^0 \int_{-\mathtt d}^0 \<r_{yy} (t,  S(t),S_t)(V(t)- \pi(t))p^{a,2}(t)y^{s,h,S}(t+\theta), y^{s,h,S}(t+\eta)\>\,\mu_r(d\theta)\,\mu_r(d\eta)\,dt 
      \\\nonumber &- \E ^{\F_s} \int_s^T \int_{-\mathtt d}^0 \int_{-\mathtt d}^0\<\pi(t)b_{yy} (t,  S(t),S_t)p^{a,2}(t)y^{s,h,S}(t+\theta), y^{s,h,S}(t+\eta)\>\,\mu_b(d\theta)\,\mu_b(d\eta)\,dt \\\nonumber &- \E ^{\F_s} \int_s^T \int_{-\mathtt d}^0 \int_{-\mathtt d}^0\<\pi(t)\sigma_{yy} (t,  S(t),S_t)q^{a,2}(t)y^{s,h,S}(t+\theta), y^{s,h,S}(t+\eta)\>\,\mu_\sigma(d\theta)\,\mu_\sigma(d\eta)\,dt 
      \\ \nonumber &
      - \E ^{\F_s} \int_s^T \<(U_1)_{vv}y^{s,h,V}(t), y^{s,h,V}(t)\>\,dt -2  \E ^{\F_s} \int_s^T \<r_{x} (t,  S(t),S_t)p^{a,2}(t)y^{s,h,S}(t), y^{s,h,V}(t)\>\,dt \\ \nonumber &
      -2  \E ^{\F_s} \int_s^T \int_{-\mathtt d}^0 \<r_{y} (t,  S(t),S_t)p^{a,2}(t)y^{s,h,S}(t+\theta), y^{s,h,V}(t)\>\,\mu_r(d\theta)dt \\ \nonumber\end{align}
      \begin{align*}
       &
      -2  \E ^{\F_s} \int_s^T\int_{-\mathtt d}^0  \<r_{xy} (t,  S(t),S_t)(V(t)- \pi(t))p^{a,2}(t)y^{s,h,S}(t+\theta), y^{s,h,S}(t)\>\,\mu_r(d\theta) dt 
      \\ \nonumber &-2  \E ^{\F_s}  \int_s^T \int_{-\mathtt d}^0\<\pi(t)b_{xy} (t,  S(t),S_t)p^{a,2}(t)y^{s,h,S}(t+\theta), y^{s,h,S}(t)\>\,\mu_b(d\theta)dt \\ \nonumber &
      -2  \E ^{\F_s}  \int_s^T \int_{-\mathtt d}^0\<\pi(t)\sigma_{xy} (t,  S(t),S_t)q^{a,2}(t)y^{s,h,S}(t+\theta), y^{s,h,S}(t)\>\,\mu_\sigma(d\theta)dt 
    \end{align*}
From the maximum principle stated in Theorem \ref{teoMaxPrinc-gen} we deduce the following condition on the optimal strategy for the present problem: notice that the optimality condition can be given only in terms of the pair of processes $(p^{a,2},q^{a,2})$, being $(p^{a,1},q^{a,1})$ identically equal to $0$.
\begin{theorem}\label{maxprinc-findim_executiondelay}
 Let Hypothesis \ref{ip: b_r_sigma} holds true. Let $(p^{a,2},q^{a,2})$ be the unique solution of the ABSDE (\ref{ABSDE-aggiuntaI_execution delay2})  and $P_{0,0}$ given in \eqref{teoP00formula-example}.
Then the following variational  inequality holds true $\forall\, (\pi, c)\in \R\times \mathcal C$:\vspace{-0.2 truecm}
\begin{multline}
 \calh^a(t, x(t), \pi(t),c(t),p^a(t),q^a(t))-  \calh^a(t,x(t), \pi,c,p^a(t),q^a(t))
\\ \nonumber +\< P_{0,0}(t) (\pi(t)\sigma (t,   S(t), S_t)- \pi\sigma (t,S(t),S_t)),\pi(t)\sigma (t, S(t),S_t)- \pi\sigma (t,S(t),S_t)\>\leq 0
\end{multline}
$dt\times \P\; a.s.,$ where $\calh^a $ is defined in \eqref{hamiltonian-example}.
\end{theorem}

\section{Proof of propositions \ref{prop:Ito}, \ref{prop.conv.yn-zn} and \ref{prop.conv yhn}.}\label{appendix}

{\bf Proof of Proposition \ref{prop:Ito}} By applying the Ito formula we compute
\begin{align*}
 &d( p^a(t) y^\e(t))\\ \nonumber
 &=d (p^a(t)) y^\e(t)+p^a(t)\,d y^\e(t)-q(t)(\sigma_x(t) y^\e(t) + \sigma_y(t) \int_{-\mathtt{d}}^0 y^ \e (t+ \theta) \mu_{\sigma}(d\theta) +\delta\sigma(t))dt  \\ \nonumber
 &=-(\ell_x(t)+  \dis \int_{-\mathtt{d}}^0\E^{\Fcal_t}\ell_y(t-\theta)\mu_\ell(d\theta) +\dis b_x(t) {p}^a(t)+ \dis \int_{-\mathtt{d}}^0 \, \E^{\Fcal_t}b_y(t-\theta)p^a(t-\theta)\mu_b(d\theta)  \\& +\sigma_x(t) {q}^a(t) ) y^\e(t)\,dt 
-\int_{-\mathtt{d}}^0 \E^{\Fcal_t}\sigma_y(t-\theta) {q}^a(t-\theta)\mu_\sigma(d\theta) \, d\theta \,y^\e(t)\,dt +q^a(t), y^\e(t)\, dW(t)
 \>\\ \nonumber
 &
 +p^a(t)( b_x(t) y^\e(t) + b_y(t) \int_{-\mathtt{d}} ^ 0 y^ \e (t+ \theta) {\mu} _b(d\theta) )\,dt\\ \nonumber
 &-q^a(t)(\sigma_x(t) y^\e(t) + \sigma_y(t) \int_{-\mathtt{d}}^0 y^ \e (t+ \theta) \mu_{\sigma}(d\theta) +\delta\sigma(t))dt \hspace{5truecm}
\end{align*}
Integrating between $0$ and $T$ and taking expectation we obtain\vspace{-0.2truecm}
\begin{equation*}
\E h_x(x(T)), y^\e(T)=\E\dis\int_0^T[\ell_x(s)y^\e(s)+q^a(s)\delta\sigma(s) +\dis\int_{-\mathtt{d}}^0 \ell_y(s)   y^\e(s+\theta)  \mu_\ell(d\theta) ]\,ds
\end{equation*}
since we have that\vspace{-0.2truecm}
\begin{align}\label{cancel-b}
    &-\E\dis\int_0^T\dis \int_{-\mathtt{d}}^0 \, \E^{\Fcal_t}b_y(t-\theta)p^a(t-\theta)\mu_b(d\theta) y^\e(t)\,dt+\E\dis\int_0^T\dis\int_{-\mathtt{d}}^0  p^a(t)b_y(t) y^ \e (t+ \theta) {\mu} _b(d\theta) \,dt\\
     &=-\E\dis\int_0^T\dis \int_{-\mathtt{d}}^0 \, \E^{\Fcal_t}b_y(\tau)p^a(\tau)\mu_b(d\theta) y^\e(\tau+\theta)\,d\tau+\E\dis\int_0^T\dis\int_{-\mathtt{d}}^0  p^a(t)b_y(t) y^ \e (t+ \theta) {\mu} _b(d\theta)\, dt=0 \nonumber
\end{align}
where in the first integral we have made the change of time $\tau=t-\theta$, taking into account that $y^\e(\tau+\theta)\equiv 0$ if $ 0 \leq \tau\leq -\theta$ and $p^a(t-\theta)\equiv 0$ if $t>T+\theta$; analogously we have
\begin{align}\label{cancel-sigma}
    &-\E\dis\int_0^T\dis \int_{-\mathtt{d}}^0 \, \E^{\Fcal_t}\sigma_y(t-\theta)q^a(t-\theta)\mu_\sigma(d\theta)y^\e(t)\,dt+\E\dis\int_0^T\dis\int_{-\mathtt{d}}^0  q^a(t)\sigma_y(t) y^ \e (t+ \theta) {\mu} _\sigma(d\theta)\, dt\\
    & =-\E\dis\int_0^T\dis \int_{-\mathtt{d}}^0 \, \E^{\Fcal_t}\sigma_y(\tau)q^a(\tau)\mu_\sigma(d\theta) y^\e(\tau+\theta)\,d\tau+\E\dis\int_0^T\dis\int_{-\mathtt{d}}^0  q^a(t)\sigma_y(t) y^ \e (t+ \theta) {\mu} _\sigma(d\theta)\, dt=0 \nonumber
\end{align}
Next by applying the Ito formula we compute 
\begin{align*}
 d&( p^a(t) z^\e(t))=(d p^a(t)) z^\e(t)+p^a(t)(d z^\e(t))-q^a(t)(\sigma_x(t) z^\e(t) \\ \nonumber
 &\,+\sigma_y(t) \int_{-\mathtt{d}}^0 z^ \e (t+ \theta) \mu_{\sigma}(d\theta) \, dt+ \delta {\sigma}_x(t) y^\e(t) )\,dt\\ \nonumber
 &\,-q^a(t)( \displaystyle  \delta {\sigma}_y(t) \displaystyle \int_{-\mathtt{d}}^ 0y^ {\e} (t+ \theta) \mu_\sigma(d\theta) -\int_{-\mathtt{d}}^0\sigma_{xy}(t)y^\e(t)  y^\e(t+\theta){\mu}_\sigma(d\theta))
 \, dt \\ &  -\frac{1}{2}q^a(t)(\sigma_{xx}(t)y^\e(t) +\int_{-\mathtt{d}}^ 0\int_{-\mathtt{d}}^ 0\sigma_{yy}(t) y^\e(t+\theta_1)y^\e(t+\theta_2) {\mu} _\sigma(d\theta_1) {\mu} _\sigma(d\theta_2)) \,dt
    \end{align*}
Integrating between $0$ and $T$, taking expectation and arguing as in \eqref{cancel-b} and \eqref{cancel-sigma} we obtain
\begin{align*}
&\E h_x(x(T)) z^\e(T)=\E\dis\int_0^T[\ell_x(s)z^\e(s) +
\int_{-\mathtt{d}}^0 \ell_y(s)   z^\e(s+\theta)  \mu_\ell(d\theta)
] \, ds\\&+ \E\dis\int_0^Tq^a(s)[] \delta {\sigma}_x(s) y^\e(s)+   \displaystyle  \delta {\sigma}_y(s) \displaystyle \int_{-\mathtt{d}}^ 0y^ {\e} (s+ \theta) \mu_\sigma(d\theta)]\,ds\\
&\quad+ \E\dis\int_0^T\frac{1}{2}p^a(s)[ b_{xx}(s)y^\e(s)y^\e(s)+\int_{-\mathtt{d}}^ 0\int_{-\mathtt{d}}^ 0b_{yy}(s) y^\e(s+\theta_1)y^\e(t+\theta_2){\mu} _b(d\theta_1) {\mu} _b(d\theta_2)]\,ds\\
&\quad +\E\dis\int_0^Tp^a(s)[\int_{-\mathtt{d}}^0 b_{xy}(s)y^\e(s)  y^{\e}(s+\theta)\mu_b(d\theta)+\delta b(s)]\,ds\\
 &\quad+ \E\dis\int_0^T\frac{1}{2}q^a(s)[ \sigma_{xx}(s)y^\e(s)y^\e(s)+\int_{-\mathtt{d}}^ 0\int_{-\mathtt{d}}^ 0\sigma_{yy}(s) y^\e(s+\theta_1)y^\e(s+\theta_2){\mu} _\sigma(d\theta_1) {\mu} _\sigma(d\theta_2)]\,ds\\
&\quad +\E\dis\int_0^Tq^a(s)[\int_{-\mathtt{d}}^0 \sigma_{xy}(s)y^\e(s)  y^{\e}(s+\theta)\mu_\sigma(d\theta)]\,ds
\end{align*}
and this concludes the proof. \qed

\medskip

{\bf Proof of Proposition \ref{prop.conv.yn-zn} }
Existence  and regularity of the solutions  follows from \cite{Moha}.
Let us move to the convergence result.
Comparing the equations \eqref{agg.prima.fin.dim}  and \eqref{yn}, by standard estimates on the difference of the two equations, we get that:
\begin{align}\label{eq.stima.diff}
\E &\sup_{-\mathtt{d} \leq t \leq r} |y^{n,\e}(t)- y^{\e}(t)|^2 \leq C \E \int_0^r \Big|\int_{-\mathtt{d}}^0 y^{n,\e}(s+ \theta ) f^n_b(\theta) \, d\theta - \int_{-\mathtt{d}}^0 y^{n,\e}(s+ \theta ) \mu_b(d\theta) \Big|^ 2 ds\\ \nonumber& \hspace{4cm}+C \E \int_0^r \Big|\int_{-\mathtt{d}}^0 y^\e(s+ \theta ) f^n_\sigma(\theta) \, d\theta - \int_{-\mathtt{d}}^0 y^\e(s+ \theta ) \mu_\sigma(d\theta) \Big|^ 2 ds\\ \nonumber
&\leq C \E \int_{-\mathtt{d}}^r  \sup_{-\mathtt{d} \leq s \leq \sigma} \Big|y^{n,\e}(s)- y^{\e}(s)\Big|^2 \, ds \nonumber  +  
C \E \int_0^T\Big|\int_{-\mathtt{d}}^0 y^\e(s+ \theta ) f^n_b(\theta) \, d\theta - \int_{-\mathtt{d}}^0 y^\e(s+ \theta ) \mu_b(d\theta) \Big|^ 2 ds\\ \nonumber&\hspace{4cm} +C \E \int_0^T \Big|\int_{-\mathtt{d}}^0 y^\e(s+ \theta ) f^n_\sigma(\theta) \, d\theta - \int_{-\mathtt{d}}^0 y^\e(s+ \theta ) \mu_\sigma(d\theta) \Big|^ 2 ds 
\end{align}
To get the second inequality we have added and subtracted the  terms $\displaystyle\int_{-\mathtt{d}}^0 y^{\e} (t+ \theta) f^n _b (\theta)d\theta $ and $\displaystyle\int_{-\mathtt{d}}^0 y^{\e} (t+ \theta) f^n _\sigma (\theta)d\theta $ and we have taken into account that
that $|\mu^n| \leq K$ for some $K >0$ independent of $n$, so that
\begin{align*} 
 \E \int_0^r \Big|\int_{-\mathtt{d}}^0 (y^{n,\e}(s+ \theta )-y^{\e}(s+ \theta )) &f^n_b(\theta) \, d\theta\Big|^2\,ds+\E \int_0^r\Big|\int_{-\mathtt{d}}^0 (y^{n,\e}(s+ \theta )-y^{\e}(s+ \theta )) f^n_\sigma(\theta) \, d\theta\Big|^2\,ds\\
&\leq C \E \int_{-\mathtt{d}}^r  \sup_{-\mathtt{d} \leq s \leq \sigma} \Big|y^{n,\e}(s)- y^{\e}(s)\Big|^2 
\end{align*}
Taking into account we have that $\vert\mu\vert\leq K$ ans also $\vert\mu_n\vert \leq K$, uniformly with respect to $n$, by the dominated convergence theorem, and noting that $y^\e(\cdot)$ is continuous, the last integrals in \eqref{eq.stima.diff}  go to $0$ as $n\rightarrow +\infty$. So by applying the Gromwall Lemma in \eqref{eq.stima.diff} we deduce convergence of $y^{n,\e}$ \eqref{conv y}.
\newline Once we have this convergence result we study the difference $z^\e(t)-z^{n,\e}(t)$. As before we get: 

\begin{align*}
&\E \sup_{-\mathtt{d}\leq t \leq r} |z^{n,\e}(t)- z^{\e}(t)|^2 \leq C \E \int_{-d}^r  \sup_{-\mathtt{d} \leq s \leq \sigma} \Big|z^{n,\e}(s)- z^{\e}(s)\Big|^2 \, ds  \\
&+c \E \int_0^T\left[ \Big|\int_{-\mathtt{d}}^0 z^\e(s+ \theta ) f^n_b(\theta) \, d\theta - \int_{-\mathtt{d}}^0 z^\e(s+ \theta ) \mu_b(d\theta) \Big|^ 2  + 
        \Big|\int_{-\mathtt{d}}^0 z^\e(s+ \theta ) f^n_\sigma(\theta) \, d\theta - \int_{-\mathtt{d}}^0 z^\e(s+ \theta ) \mu_\sigma(d\theta) \Big|^ 2 \, ds\right]
    \\ &  + c \E \int_0^T  \Big| \delta {\sigma}_y(s) \Big[ \int_{-\mathtt{d}}^ 0y^ {n,\e} (s+ \theta) f^n_\sigma(\theta)d\theta -  \int_{-\mathtt{d}}^ 0y^ {n,\e} (s+ \theta) \mu_\sigma(d\theta) \Big]  \Big|^2 ds   \\
    &+c\E\displaystyle\int_0^T\vert y^{n,\e}(s)y^{n,\e}(s)-y^{\e}(s)y^{\e}(s)\vert^2\,ds\\
    &+c\E\displaystyle\int_0^T\vert\int_{-\mathtt{d}}^ 0 y^{n,\e}(s) y^{n,\e}(s+\theta)f^n_b(\theta)d\theta-\int_{-d} ^ 0y^{\e}(s) y^{\e}(s+\theta)\mu_b(\theta)\vert^2\,dt
    \\   &  +c \E \int_0^T \Big| \int_{-\mathtt{d}}^0  \int_{-\mathtt{d}}^0 
y^{n,\e}(s+ \theta_1 )y^{n,\e}(s+ \theta_2 ) f^n_b(\theta_1)f^n_b(\theta_2) \, d\theta _1  d\theta _2-  \int_{-\mathtt{d}}^0  \int_{-\mathtt{d}}^0 
y^\e(s+ \theta_1 )y^\e(s+ \theta_2 )  \mu_b(d\theta_1) \mu_b(d\theta_2) \Big| ^2 \, ds \\ 
    &+c\E\displaystyle\int_0^T\vert\int_{-\mathtt{d}}^ 0 y^{n,\e}(s) y^{n,\e}(s+\theta)f^n_\sigma(\theta)d\theta-\int_{-d} ^ 0y^{\e}(s) y^{\e}(s+\theta)\mu_\sigma(\theta)\vert^2\,ds +\\
    & c \E \int_0^T \Big| \int_{-\mathtt{d}}^0  \int_{-\mathtt{d}}^0 y^{n,\e}(s+ \theta_1 )y^{n,\e}(s+ \theta_2 ) f^n_\sigma(\theta_1)f^n_\sigma(\theta_2) \, d\theta _1  d\theta _2-  \int_{-\mathtt{d}}^0  \int_{-\mathtt{d}}^0 y^\e(s+ \theta_1 )y^\e(s+ \theta_2 )  \mu_\sigma(d\theta_1) \mu_\sigma(d\theta_2) \Big| ^2 \, ds
\end{align*}
Let us see how to treat the last term, since the other terms can be handled similarly:
\begin{align*}
 & 
\E \int_0^T \Big| \int_{-\mathtt{d}}^0  \int_{-\mathtt{d}}^0 
y^{n,\e}(s+ \theta_1 )y^{n,\e}(s+ \theta_2 ) f^n_\sigma(\theta_1)f^n_\sigma(\theta_2) \, d\theta _1  d\theta _2-  \int_{-\mathtt{d}}^0  \int_{-\mathtt{d}}^0 
y^\e(s+ \theta_1 )y^\e(s+ \theta_2 )  \mu_\sigma(d\theta_1) \mu_\sigma(d\theta_2) \Big| ^2 \, ds  \\
& 
\leq  c \Big [\E \int_0^T \Big| \int_{-\mathtt{d}}^0
y^{n,\e}(s+ \theta_1 ) f^n_\sigma(\theta_1) d\theta_1\int_{-\mathtt{d}}^0 \Big( y^{n,\e}(s+ \theta_2 ) -  
y^\e(s+ \theta_2 ) \Big) f^n_\sigma(\theta_2)  d\theta _2 \Big| ^2 \, ds  \\  &  +\E \int_0^T \Big| \int_{-\mathtt{d}}^0
y^{n,\e}(s+ \theta_1) f^n_\sigma(\theta_1) d\theta_1\Big( \int_{-\mathtt{d}}^0 y^\e(s+ \theta_2 ) f^n_\sigma(\theta_2)  d\theta _2- 
y^\e(s+ \theta_2 )  \mu_\sigma(d\theta_2)  \Big)\Big| ^2 \, ds \\&  +\E \int_0^T \Big| \int_{-\mathtt{d}}^0
y^\e(s+ \theta_2 ) \mu_\sigma(\theta_2) d\theta_2 \int_{-\mathtt{d}}^0  \Big(y^{n,\e}(s+ \theta_1 ) - 
y^\e(s+ \theta_1 )  \Big) f^n_\sigma(\theta_1)  d\theta _1\Big| ^2 \, ds \\&   
+\E \int_0^T \Big| \int_{-\mathtt{d}}^0
y^\e(s+ \theta_2 ) \mu_\sigma(\theta_2) d\theta_2\Big( \int_{-\mathtt{d}}^0 y^\e(s+ \theta_1 ) f^n_\sigma(\theta_1)  d\theta _1- 
y^\e(s+ \theta_1 )  \mu_\sigma(d\theta_1)  \Big)\Big| ^2 \, ds \Big ]
\end{align*}
Notice that 

\begin{align*}
& \E \int_0^T \Big| \int_{-\mathtt{d}}^0
y^\e(s+ \theta_2 ) \mu_\sigma(\theta_2) d\theta_2 \int_{-\mathtt{d}}^0  \Big(y^{n,\e}(s+ \theta_1 ) - 
y^\e(s+ \theta_1 )  \Big) f^n_\sigma(\theta_1)  d\theta _1\Big| ^2 \, ds \\&  \leq  
\E \Big [ \Big(\int_0^T \Big| \int_{-\mathtt{d}}^0
y^\e(s+ \theta_2 ) \mu_\sigma(\theta_2) d\theta_2 \Big| ^2  \, ds \Big)^ {1/2} \Big(\int_0^T \Big| \int_{-\mathtt{d}}^0  \Big(y^{n,\e}(s+ \theta_1 ) - 
y^\e(s+ \theta_1 )  \Big) f^n_\sigma(\theta_1)  d\theta _1\Big| ^2 \, ds  \Big)^ {1/2} \Big]
\\& \leq  c \E \Big [\sup_{t \in [-\mathtt{d},T]} | y^{\e}(t)|   \Big(\int_0^T\Big| \int_{-\mathtt{d}}^0  \Big(y^{n,\e}(s+ \theta_1 ) - 
y^\e(s+ \theta_1 )  \Big) f^n_\sigma(\theta_1)  d\theta _1\Big|^2 \, ds  \Big)^ {1/2} \Big]
\\& \leq c \,  [\E\sup_{t \in [-\mathtt{d},T]} | y^{\e}(t)|^2 ]^{1/2} [\E\sup_{t \in [-\mathtt{d},T]} | y^{n,\e}(t) - 
y^\e(t)|^2 ]^{1/2}
\end{align*}
where $c$ may change from line to line but will never  depend on $n$ nor on $\e$.
Moreover\vspace{-0.2truecm}
\begin{align*}
 &  \Big| \int_{-\mathtt{d}}^0
y^{n,\e}(s+ \theta_1 ) f^n_\sigma(\theta_1) d\theta_1\Big( \int_{-\mathtt{d}}^0 y^\e(s+ \theta_2 ) f^n_\sigma(\theta_2)  d\theta _2- 
y^\e(s+ \theta_2 )  \mu_\sigma(d\theta_2)  \Big)\Big| ^2   \\&
\leq  C(K) \Big [\sup_{-\mathtt{d} \leq t \leq T} |y^{\e}(t)|^4+ \sup_{-\mathtt{d} \leq t \leq T} |y^{n,\e}(t)|^4 \Big]
\end{align*}
where $K$ is given in lemma \ref{lemma:approx-measure}.
Hence the convergence of $z^{n,\e}$ in \eqref{conv y} follows from \eqref{stima sup yn}, the Gronwall Lemma and the Dominated Convergence Theorem, as for the convergence of $y^{n,\e}$.
\qed

\medskip
{\bf Proof of Proposition \ref{prop.conv yhn}.}
As in the previous proof, estimate  \eqref{eq.stima.diff}, we arrive at considering\vspace{-0.2truecm}
\begin{align}\label{eq.stima.diff II}
\E \sup_{-\mathtt{d} \leq t \leq r} |y^{n,s,h}(t)- y^{s,h}(t)|^2 & \leq C 
\E \int_{-\mathtt{d}}^r  \sup_{-\mathtt{d} \leq \tau \leq \sigma} \Big|y^{n,s,h}(\tau)- y^{s,h}(\tau)\Big|^2 \, d\sigma \nonumber \\& +  
C \E \int_0^T \Big|\int_{-\mathtt{d}}^0 y^{s,h}(\tau+ \theta ) f^n_b(\theta) \, d\theta - \int_{-\mathtt{d}}^0 y^{s,h}(\tau+ \theta ) \mu_b(d\theta) \Big|^ 2 ds\\ \nonumber& +C \E \int_0^T \Big|\int_{-\mathtt{d}}^0 y^{s,h}(\tau+ \theta ) f^n_\sigma(\theta) \, d\theta - \int_{-\mathtt{d}}^0 y^{s,h}(\tau+ \theta ) \mu_\sigma(d\theta) \Big|^ 2 d\tau .
\end{align}
To get these estimates we have added and subtracted the  terms $\E\Big|\displaystyle\int_0^T\displaystyle\int_{-\mathtt{d}}^0 y^{s,h}(\tau+ \theta) f^n _b (\theta)d\theta d\tau \Big|^2$ and $\E\Big|\displaystyle\int_0^T\displaystyle\int_{-\mathtt{d}}^0 y^{s,h}(\tau+ \theta) f^n _\sigma (\theta)d\theta d\tau \Big|^2$ and we have taken into account that
that $|\mu^n| \leq K$ for some $K >0$ independent of $n$, so that\vspace{-0.2truecm}
\begin{align*} 
 \E \int_0^T \Big|\int_{-\mathtt{d}}^0 (y^{n,s,h}(\tau+ \theta )&-y^{s,h}(\tau+ \theta )) f^n_b(\theta) \, d\theta\Big|^2 d\tau+\E \int_0^T \Big|\int_{-\mathtt{d}}^0 (y^{n,\e}(\tau+ \theta )-y^{\e}(\tau+ \theta )) f^n_\sigma(\theta) \, d\theta\Big|^2d\tau\\
&\leq C \E \int_{-\mathtt{d}}^r  \sup_{-\mathtt{d} \leq \tau \leq \sigma} \Big|y^{n,s,h}(\tau)- y^{s,h}(\tau)\Big|^2 d\tau
\end{align*}
Unlike in the previous proof, the trajectories of $y^{s,h}$ have a discontinuity in $t=s$, so in order to exploit the weak convergence argument, we introduce the following  localizing sequence of smooth functions $\phi_s^{\delta}$, $ \delta >0$, such that $0\leq \phi_s^{\delta} \leq 1$ and\vspace{-0.2truecm}
\begin{equation}
    \phi_s^{\delta}(t)= \left \{  \begin{array}{ll}
       1   & t \in [s+ 2 \delta, T],  \\
        0 &  t \in [-\mathtt{d}, s+ \delta]
    \end{array} \right .
\end{equation}
Now we set $y^{\delta,s,h} (t):= y^{s,h} (t) \phi_s^{\delta}(t)$, $ t \in [-\mathtt{d},T]$, now
$y^{\delta,s,h} (t) \in C([-\mathtt{d},T])$. 

Thus, let see the third term at the L.H.S. of \eqref{eq.stima.diff II} (the second term at L.H.S. of \ref{eq.stima.diff II} can be treated in the same way)\vspace{-0.2truecm}
\begin{align}\label{A2-third}
& \E \int_0^T \Big|\int_{-\mathtt{d}}^0 y^{s,h}(t+ \theta ) f^n_\sigma(\theta) \, d\theta - \int_{-\mathtt{d}}^0 y^{s,h}(t+ \theta ) \mu_\sigma(d\theta) \Big|^ 2 dt  \\ \nonumber&   \leq  
 c \Big \{  \E \int_0^T \Big|\int_{-\mathtt{d}}^0[ y^{s,h}(t+ \theta)- y^{\delta, s,h}(t+ \theta ) ] f^n_\sigma(\theta) \, d\theta \Big|^ 2 dt    +
\\ \nonumber&  \E \int_0^T \Big|\int_{-\mathtt{d}}^0 y^{\delta,s,h}(t+ \theta ) f^n_b(\theta) \, d\theta - \int_{-\mathtt{d}}^0 y^{\delta,s,h}(t+ \theta ) \mu_\sigma(d\theta) \Big|^ 2 dt  +
\\ \nonumber& \E \int_0^T \Big|\int_{-\mathtt{d}}^0[ y^{s,h}(t+ \theta)- y^{\delta,s,h}(t+ \theta ) ] \mu_\sigma(d \theta)\Big|^ 2 dt  \Big\}
\end{align}
We estimate the first term of \eqref{A2-third}, \vspace{-0.2truecm}
\begin{align}\label{B1}
  & \E \int_0^T \Big|\int_{-\mathtt{d}}^0[ y^{s,h}(t+ \theta)- y^{\delta, s,h}(t+ \theta ) ] f^n_\sigma(\theta) \, d\theta \Big|^ 2 dt    \\  \nonumber & \leq   \E \int_0^T \Big|\int_{-\mathtt{d}}^0[ y^{s,h}(t+ \theta)I_{[s,s+ 2\delta]}(t+ \theta ) ] f^n_\sigma(\theta) \, d\theta \Big|^ 2 dt
   \\ \nonumber& =   \E \int_0^T  \Big|\int_{-\mathtt{d}}^0[ y^{s,h}(t+ \theta_1)I_{[s,s+ 2\delta]}(t+ \theta_1) ] f^n_\sigma(\theta) \, d\theta_1 \Big| \Big|\int_{-\mathtt{d}}^0[ y^{s,h}(t+ \theta)I_{[s,s+ 2\delta]}(t+ \theta ) ] f^n_\sigma(\theta) \, d\theta \Big|  dt 
  \\ \nonumber
& \leq  \E\sup_{-\mathtt{d}\leq \rho \leq T}  | y^{s,h}(\rho)|^ 2 \mu^n([-\mathtt{d},0])  \int_0^T  \int_{-\mathtt{d}}^0 I_{[s,s+ 2\delta]} (t+ \theta)  f^n_\sigma(\theta) \, d\theta \, dt 
\\ \nonumber& =
|\mu^n|([-\mathtt{d},0]) \E\sup_{-\mathtt{d}\leq \rho \leq T}  | y^{s,h}(\rho)| ^ 2  \int_{-\mathtt{d}}^0 (\int_{(s-\theta) \vee 0}^{ (s + \delta -\theta) \wedge T }\, dt)   \  f^n_\sigma(\theta) d\theta 
\\ \nonumber& \leq  \delta
|\mu^n|([-\mathtt{d},0])^2  \E\sup_{-\mathtt{d}\leq \rho \leq T}  | y^{s,h}(\rho)| ^ 2 \leq C \delta \E\sup_{-\mathtt{d}\leq \rho \leq T}  | y^{s,h}(\rho)| ^ 2 
\end{align}
In a similar way we can get
\begin{equation}\label{B2}
    \E \int_0^T \Big|\int_{-\mathtt{d}}^0[ y^{s,h}(t+ \theta)- y^{\delta,s,h}(t+ \theta ) ] \mu_\sigma(d \theta)\Big|^ 2 dt \leq C \delta \E\sup_{-\mathtt{d}\leq \rho \leq T}  | y^{s,h}(\rho)| ^ 2 
\end{equation}
Using all these estimates in \eqref{eq.stima.diff II} and the Gromwall Lemma as well, we get that:
\begin{align*}
&\E \sup_{-\mathtt{d} \leq t \leq r} |y^{n,s,h}(t)- y^{s,h}(t)|^2 \leq  \\ &
 C \Big[ \delta \,  \E\sup_{-\mathtt{d}\leq \rho \leq T}  | y^{s,h}(\rho)| ^ 2  +
  \E \int_0^T \Big|\int_{-\mathtt{d}}^0 y^{\delta,s,h}(t+ \theta ) f^n_b(\theta) \, d\theta - \int_{-\mathtt{d}}^0 y^{\delta,s,h}(t+ \theta ) \mu_b(d\theta) \Big|^ 2 dt   \Big],
\end{align*}
then, by the weak convergence of $\mu^n$ to $\mu$ and since $\delta$ can be chosen as small as we want,  we deduce that
\begin{align}\label{limite debole delta}
& \lim_{n \to \infty}\E \sup_{-\mathtt{d} \leq t \leq r} |y^{n,s,h}(t)- y^{s,h}(t)|^2 \\& \nonumber\leq 
  \lim_{n \to \infty}\E \int_0^T \Big|\int_{-\mathtt{d}}^0 y^{\delta,s,h}(t+ \theta ) f^n_b(\theta) \, d\theta - \int_{-\mathtt{d}}^0 y^{\delta,s,h}(t+ \theta ) \mu_b(d\theta) \Big|^ 2 dt +\bar C \delta =0.
\end{align}

Eventually in view of 
 \eqref{B1}, \eqref{B2} and \eqref{limite debole delta} we get:
\begin{align}
  \lim_{n \to \infty}\E \int_0^T \Big|\int_{-\mathtt{d}}^0 y^{s,h}(t+ \theta ) f^n_b(\theta) \, d\theta - \int_{-\mathtt{d}}^0 y^{s,h}(t+ \theta ) \mu_b(d\theta) \Big| =0.
\end{align}
\qed

\end{document}